\documentclass[12pt,leqno]{article}
\usepackage{amsmath,amssymb,amsthm,amscd}
\numberwithin{equation}{section} \allowdisplaybreaks
\newtheorem{theorem}{Theorem}[section]
\newtheorem{defin}{Definition}[section]
\newtheorem{prop}{Proposition}[section]
\newtheorem{corol}{Corollary}[section]
\newtheorem{lemma}{Lemma}[section]

\newtheorem{rem}{Remark}[section]
\newtheorem{example}{Example}[section]
\begin{document}
\title{Reduction and submanifolds of generalized complex manifolds}
\author{Izu Vaisman\footnote{vaisman@math.haifa.ac.il}\\
{\small Department of Mathematics, University of Haifa, Israel}}
\date{}
\maketitle
{\def\thefootnote{*}\footnotetext[1]%
{{\it 2000 Mathematics Subject Classification: 53D17, 53C56}.
}}
\begin{center} \begin{minipage}{12cm}
A{\footnotesize BSTRACT. We recall the presentation of the
generalized, complex structures by classical tensor fields, while
noticing that one has a similar presentation and the same
integrability conditions for generalized, paracomplex and subtangent
structures. This presentation shows that the generalized, complex,
paracomplex and subtangent structures belong to the realm of Poisson
geometry. Then, we prove geometric reduction theorems of
Marsden-Ratiu and Marsden-Weinstein type for the mentioned
generalized structures and give the characterization of the
submanifolds that inherit an induced structure via the corresponding
classical tensor fields.}
\end{minipage}
\end{center}
\vspace{5mm}
The study of generalized, complex structures is a recent subject
that was started by N. J. Hitchin \cite{Ht1} and M. Gualtieri
\cite{Galt} and was continued by several authors
\cite{{AB},{BB},{Cr}, {Ht2},{IW},{LMTZ},{Wade}}. The subject is
motivated by the fact that generalized complex manifolds appear as
target manifolds of $\sigma$-models with supersymmetries
\cite{LMTZ}.

The framework of the present paper is the $C^\infty$-category and
$M$ is a differentiable manifold. The generalized complex structures
are defined like the classical complex structures but, with the
tangent bundle $TM$ replaced by $T^{big}M=TM\oplus T^*M$ and with
the Lie bracket of vector fields replaced by the Courant bracket
\cite{{CrW},{C}}. Lindstr\"om-Minasian-Tomasiello-Zabzine
\cite{LMTZ} and Crainic \cite{Cr} gave the full interpretation of a
generalized, complex structure by means of classical tensor fields.
Particularly, among these tensor fields there is a Poisson bivector
field, which was also discovered in \cite{Galt} and \cite{AB}.

The fact that a generalized, complex manifold has an underlying
Poisson structure justifies the study of the Poisson geometry of the
generalized, complex manifolds. In particular, in \cite{Cr} the
subject is integrability of a generalized, complex manifold to a
certain type of symplectic groupoid. In the present paper we will
discuss reduction and submanifolds from the Poisson point of view.
In brief, the content of the paper is as follows.

In section 1, we present the basics of the theory of the
generalized, complex structures using the corresponding classical
tensor fields. The content of this section is not original, and its
length is justified by the fact that the whole theory is pretty new
and, presumably, not very popular yet. However, the section also
contains some novelties: the results are formulated for three kinds
of generalized structures, complex, paracomplex and subtangent; we
show the connection between the generalized structures and
Poisson-Nijenhuis structures\footnote{Recently, I learned from P. Xu
that he has also indicated such a connection in his lecture at the
conference in Trieste, Italy, July 2005.}; we indicate the algebraic
expression of a generalized structure along a symplectic leaf of its
Poisson structure; finally, we refer to the possible Lie groups with
a compatible generalized structure.

In Section 2 we discuss reduction of generalized structures. We
start with the Marsden-Ratiu definition of the geometric reduction
of a Poisson structure via a submanifold and a control vector bundle
\cite{MRt} and prove a geometric reduction theorem for generalized
structures. Then, we particularize the theorem for interesting
special control bundles. In particular, we obtain a corollary which
is a Marsden-Weinstein reduction theorem for generalized structures.

Finally, in Section 3 we discuss the notion of a submanifold defined
by Ben-Bassat-Boyarchenko in \cite{BB}. These submanifolds inherit
an induced generalized structure. We give the characterization of
the submanifolds in the sense of \cite{BB} by means of the classical
tensor fields of the structure. In particular, the submanifolds
under consideration have to be Poisson-Dirac submanifolds in the
sense of \cite{CF}. The same characterization may also be used as a
good definition of Poisson-Nijenhuis submanifolds of a
Poisson-Nijenhuis manifold, such that the submanifold inherits an
induced hierarchy of Poisson structures.

While work on this paper was in progress, several papers on
reduction of generalized complex structures were posted on the web
\cite{{SHu},{LT},{SXu},{BG},{LT2}}. These papers provide various
ways of extending the notion of a Hamiltonian Lie group action with
an equivariant momentum map and the Marsden-Weinstein reduction
theorem to generalized complex structures. Instead, in the present
paper we extend the Marsden-Ratiu reduction theorem of \cite{MRt}.
\section{Generalized structures in classical terms}
We begin by recalling the {\it Courant bracket} $[\,,\,]:\Gamma
T^{big}M\times\Gamma T^{big}M\rightarrow\Gamma T^{big}M$ ($\Gamma$
denotes spaces of cross sections of vector bundles), which is
defined by \cite{{CrW},{C}}
\begin{equation}\label{crosetC}
[(X,\alpha),(Y,\beta)]=([X,Y],L_X\beta-L_Y\alpha
+\frac{1}{2}d(\alpha(Y)-\beta(X)),\end{equation} where $X,Y\in
\chi^1(N),\,\alpha,\beta\in\Omega^1(N)$ (we denote by $\chi^k(M)$
the space of $k$-vector fields and by $\Omega^k(M)$ the space of
differential $k$-forms on $M$). We also recall that $T^{big}M$ has
the neutral metric
\begin{equation}\label{gptCourant} g((X,\alpha),(Y,\beta))=
\frac{1}{2}(\alpha(Y)+\beta(X))
\end{equation}  and the non degenerate $2$-form
\begin{equation}\label{omegaptCourant} \omega((X,\alpha),(Y,\beta))=
\frac{1}{2}(\alpha(Y)-\beta(X)).
\end{equation}

A maximal, $g$-isotropic subbundle $L\subseteq T^{big}M$ is called
an {\it almost Dirac structure} of $M$. If $L$ is also closed by
Courant brackets it is called a {\it Dirac structure}.
\begin{rem}\label{Lreconstruit} {\rm The motivation for the study
of Dirac structures comes from the theory of constrained mechanical
systems as shown by the following facts \cite{{CrW},{C}}. An almost
Dirac structure $L$ defines a generalized distribution $D=pr_{TM}L$
(the projection $pr_{TM}$ is given by the direct sum structure of
$T^{big}M$) endowed with a $2$-form $\vartheta$ induced by $\omega$.
The structure $L$ may be reconstructed from the pair $(D,\vartheta)$
by the formula
$$ L= \{(X,\alpha)\,/\, X\in
D\,,\,\alpha|_{L_+}=\flat_{\vartheta}X\}. $$ Furthermore, $L$ is
Dirac iff $D$ is integrable and $\vartheta$ is closed along the
leaves; one says that $D$ is the {\it presymplectic foliation} of
$L$. Thus, a Dirac structure is equivalent with a generalized,
presymplectic foliation with a presymplectic form that is
differentiable on $M$.} \end{rem}

The definition of generalized, complex structures uses a
generalization of the Nijenhuis tensor. Namely, if
$\Phi\in\Gamma(End\,T^{big}M)$ one defines the {\it
Courant-Nijenhuis torsion} of $\Phi$ by
\begin{equation}\label{tensNij}
\mathcal{N}_\Phi(\mathcal{X},\mathcal{Y}) =
[\Phi\mathcal{X},\Phi\mathcal{Y}] -
\Phi[\mathcal{X},\Phi\mathcal{Y}] -
\Phi[\Phi\mathcal{X},\mathcal{Y}] +\Phi^2[\mathcal{X},\mathcal{Y}],
\end{equation} where $ \mathcal{X}=(X,\alpha),\mathcal{Y}=(Y,\beta)
\in\Gamma T^{big}M$ and the brackets are Courant brackets. In the
general case, the Courant-Nijenhuis torsion is not
$C^\infty(M)$-bilinear since, in view of the properties of the
Courant bracket \cite{C},
$\mathcal{N}_\Phi(\mathcal{X},f\mathcal{Y})$ $(f\in C^\infty(M))$
includes the terms
\begin{equation}\label{nonlinear}[g(\mathcal{X},\Phi\mathcal{Y}) +
g(\Phi\mathcal{X},\mathcal{Y})]\Phi(0,df) -
[g(\mathcal{X},\mathcal{Y})\Phi^2(0,df) +
g(\Phi\mathcal{X},\Phi\mathcal{Y})(0,df)].\end{equation} But, if
$\Phi$ is {\it $g$-skew-symmetric}, i.e.,
\begin{equation}\label{skewsym} g(\mathcal{X},\Phi\mathcal{Y}) +
g(\Phi\mathcal{X},\mathcal{Y})=0,\end{equation} and {\it
$\epsilon$-potent}, i.e.,
\begin{equation}\label{2-reductibil} \phi^2=\epsilon
Id,\hspace{1cm}\epsilon=\pm1,0,\end{equation} $\Phi$ also satisfies
the condition
\begin{equation}\label{skewsym2}
g(\Phi\mathcal{X},\Phi\mathcal{Y}) + \epsilon
g(\mathcal{X},\mathcal{Y})=0,\end{equation} and the terms
(\ref{nonlinear}) vanish.

If $\Phi$ satisfies (\ref{2-reductibil}) with $\epsilon=-1$ and
(\ref{skewsym}) $\Phi$ is equivalent with a decomposition of the
complexification $T^{big}_cM=T^{big}M\otimes_{ \mathbb{R}}
\mathbb{C}$ into a Whitney sum of conjugated complex, almost Dirac
structures, the $\pm i$-eigenbundles $L_{\pm}$ of $\Phi$, and $\Phi$
is called a {\it generalized, almost complex structure} of $M$. If
$\Phi$ satisfies (\ref{2-reductibil}) with $\epsilon=1$ and
(\ref{skewsym}) $\Phi$ is equivalent with a decomposition of
$T^{big}M$ into a Whitney sum of two maximally $g$-isotropic
subbundles, the $\pm1$-eigenbundles $E_{\pm}$ of $\Phi$ and $\Phi$
is called a {\it generalized, almost paracomplex structure} of $M$.
If $\Phi$ satisfies (\ref{2-reductibil}) with $\epsilon=0$ and
(\ref{skewsym}) we will say that $\Phi$ is a {\it generalized,
almost subtangent structure} and $im\,\Phi$, the $0$-eigenspaces
field $S$ of $\Phi$ is $g$-isotropic in $T^{big}M$. The name of this
kind of structures comes from the fact that a structure defined by a
tensor field $\Phi\in \Gamma\,End(TM)$ such that $\Phi^2=0$ and
\begin{equation}\label{conddetg}im\,\Phi=ker\,\Phi \end{equation}
is called an almost tangent structure on $M$. If the generalized,
almost subtangent structure $\Phi$ satisfies (\ref{conddetg}) $\Phi$
is a {\it generalized, almost tangent structure} of $M$ and
$S=im\,\Phi$ is an almost Dirac structure. In all the cases
mentioned above (i.e., $\epsilon=\pm1,0$), if $\mathcal{N}_\Phi=0$
the adjective ``almost" is dropped and $\Phi$ is said to be {\it
integrable}.

The following proposition gives an alternative characterization of
the generalized complex and paracomplex structures (but does not
provide a sufficient condition for the integrability of a
generalized, almost subtangent structure).
\begin{prop}\label{integrprinDirac} {\rm\cite{{Galt},{Wade}}} A generalized,
almost complex or almost paracomplex structure $\Phi$ is integrable
iff the eigenbundles of $\Phi$ are Dirac structures. If a
generalized, almost subtangent structure $\Phi$ is integrable
$im\,\Phi$ is closed by Courant brackets, and it is a Dirac
structure in the tangent case; if $im\,\Phi$ is Dirac one has
$\Phi\circ\mathcal{N}_\Phi=0$.\end{prop}
\begin{proof} Compute the values of the Courant-Nijenhuis torsion
$\mathcal{N}_{\Phi}$ on eigenvectors. For the tangent case, look at
(\ref{tensNij}) and (\ref{skewsym2}).
\end{proof}

As indicated by the title, we are interested in the generalized,
complex manifolds. However, at almost no extra cost, we get results
for all the structures mentioned above. Accordingly, we will use the
term {\it generalized (almost) c.p.s. structure}, where the letters
c,p,s stand for complex, paracomplex and subtangent, respectively.
\begin{rem}\label{Calgebroid} {\rm The above definitions may be
applied to vector bundles and Courant algebroids \cite{{LWX},{BG}}.
For instance \cite{LWX}, if $\pi$ is a Poisson bivector field on
$M$, $T^{big}M$ also has the Courant algebroid structure with anchor
$Id+\sharp_\pi$ and bracket
\begin{equation}\label{crosetCP} [(X,\alpha),(Y,\beta)]_\pi=
([X,Y]+L_\alpha Y-L_\beta X-\frac{1}{2}\sigma(\alpha(Y)-\beta(X)),
\end{equation} $$\{\alpha,\beta\}_\pi+L_X\beta-L_Y\alpha+\frac{1}{2}
d(\alpha(Y)-\beta(X)))$$ $$=
([X,Y]+i(\beta)L_X\pi-i(\alpha)L_Y\pi-\frac{1}{2}\sharp_\pi
d(\alpha(Y)-\beta(X)),$$
$$\{\alpha,\beta\}_\pi+L_X\beta-L_Y\alpha+\frac{1}{2}
d(\alpha(Y)-\beta(X))).$$ The notation is that of \cite{V-carte},
$\sigma$ is the Lichnerowicz-Poisson differential and
\begin{equation}\label{crosetforme} \{\alpha,\beta\}_\pi=
L_{\sharp_\pi\alpha}\beta - L_{\sharp_\pi\beta}\alpha -
d(\pi(\alpha,\beta)). \end{equation} This example is not very
interesting because the mapping \begin{equation}\label{isomorfism}
(X,\alpha)\mapsto(X+\sharp_P\alpha,\alpha)\end{equation} yields an
isomorphism from the new Courant algebroid to the classical Courant
algebroid, which sends the bracket (\ref{crosetCP}) to
(\ref{crosetC}) and commutes with the Nijenhuis torsion, therefore,
it sends generalized c.p.s. structures with respect to
(\ref{crosetCP}) to generalized c.p.s. structures with respect to
(\ref{crosetC}).}\end{rem}

We intend to use the interpretation of the generalized structures in
terms of classical tensor fields on $M$. For this purpose we
represent $\Phi$ in the following matrix form \cite{Galt}
\begin{equation}\label{matriceaPhi} \Phi\left(
\begin{array}{c}X\vspace{2mm}\\ \alpha \end{array}
\right) = \left(\begin{array}{cc} A&\sharp_\pi\vspace{2mm}\\
\flat_\sigma&B\end{array}\right) \left( \begin{array}{c}X\vspace{2mm}\\
\alpha \end{array}\right) \end{equation} where $(X,\alpha) \in
T^{big}M$ and, if we denote by $pr$ the natural projections and by
$\iota$ the natural embeddings, we have
$$A=pr_{TM}\circ\Phi\circ \iota_{TM}:TM \rightarrow TM,\;
\sharp_\pi=pr_{TM}\circ\Phi\circ \iota_{T^*M}:T^*M \rightarrow TM,$$
$$\flat_\sigma=pr_{T^*M}\circ\Phi\circ \iota_{TM}:TM \rightarrow T^*M,\;
B=pr_{T^*M}\circ\Phi\circ \iota_{T^*M}:T^*M \rightarrow T^*M.$$

With this notation, condition (\ref{skewsym}) is equivalent with the
following three facts:

i) $\sharp_\pi$ is defined by a bivector $\pi$ by
$\sharp_\pi\alpha=i(\alpha)\pi$,

ii) $\flat_\sigma$ is defined by a $2$-form $\sigma$ by
$\flat_\sigma X=i(X)\sigma$,

iii) $B=-^t\hspace{-1pt}A$, where $t$ denotes transposition, i.e.,
$B\alpha=-\alpha\circ A$,\vspace{2mm}\\ and we have
\begin{equation}\label{Philinie}
\Phi(X,\alpha)=(AX+\sharp_\pi\alpha,\flat_\sigma X -\alpha\circ A).
\end{equation}

Furthermore, condition (\ref{2-reductibil}) is equivalent to
\begin{equation}\label{2-reddezv} A^2=\epsilon Id -
\sharp_\pi\circ\flat_\sigma,\;\pi(\alpha\circ A,\beta)=\pi(\alpha,
\beta\circ A),\;\sigma(AX,Y)=\sigma(X,AY).\end{equation} If the
second, respectively the third, condition (\ref{2-reddezv}), holds,
$\pi$, respectively $\sigma$, is said to be {\it compatible} with
$A$.
\begin{rem}\label{obsdimpara} {\rm As a consequence of
(\ref{2-reddezv}), it follows that if a manifold $M$ has a
generalized, almost complex structure, the dimension of $M$ is even
\cite{{Galt},{BB}}. Indeed, (\ref{2-reddezv}) implies that
$(A|_{ker\,\flat_\sigma})^2=\epsilon Id$ Hence, for $\epsilon=-1$,
$dim(ker\,\flat_\sigma)$ is even. Since $dim(im\,\flat_\sigma)$ is
even too, $dim M$ is even. The same is true for generalized almost
tangent manifolds $M$ but, the argumentation is different. If $\Phi$
is the generalized almost tangent structure, there exist
decompositions $T^{big}M=im\,\Phi\oplus D$, where the terms are
maximal, $g$-isotropic and $\Phi|_D:D\rightarrow im\,\Phi$ is an
isomorphism. Hence, there exists a non degenerate $2$-form on $D$
given by $\varpi(Z_1,Z_2)=g(\Phi Z_1,Z_2)$ $(Z_1,Z_2\in D)$, and
$dim\,D=dim\,M$ must be even. On the other hand, on any manifold $M$
the decomposition $T^{big}M=TM\oplus T^*M$ is a generalized, almost
paracomplex structure while $M$ may also be
odd-dimensional.}\end{rem}

For $\Phi^2=\epsilon Id, \epsilon=-1$, the invariant computation of
the Nijenhuis torsion $ \mathcal{N}_\Phi$ with $\Phi$ given by
(\ref{matriceaPhi}) was done by Crainic \cite{Cr}. (The
corresponding computation in local coordinates appeared in
\cite{LMTZ}.) With minor adjustments, Crainic's computation also
holds for $\epsilon=1,0$ and the result is
\begin{theorem}\label{thCrainic} {\rm \cite{Cr}} The almost
c.p.s. structure $\Phi$ given by {\rm(\ref{matriceaPhi})} is
integrable iff the following conditions hold:

i) the bivector field $\pi$ defines a Poisson structure on $M$;

ii) the bracket $\{\alpha,\beta\}_\pi$ defined by
{\rm(\ref{crosetforme})} satisfies the condition
\begin{equation}\label{eqdinii}
\{\alpha,\beta\}_\pi\circ A = L_{\sharp_\pi\alpha}(\beta\circ A) -
L_{\sharp_\pi\beta}(\alpha\circ A) - d(\pi(\alpha\circ A,\beta));
\end{equation}

iii) the Nijenhuis tensor of $A$ satisfies the condition
\begin{equation}\label{Nijptintegrab} \mathcal{N}_A(X,Y) =
\sharp_\pi[i(Y)i(X)d\sigma];\end{equation}

iv) the {\it associated form}
\begin{equation}\label{sigmaA}\sigma_A(X,Y)=\sigma(AX,Y)\end{equation}
satisfies the condition
\begin{equation}\label{difsigmaA}
d\sigma_A(X,Y,Z)=\sum_{Cycl(X,Y,Z)}d\sigma(AX,Y,Z).\end{equation}
\end{theorem}

It is interesting to notice the following interpretation of
condition ii). A pair of tensor fields
$\pi\in\chi^2(M),A\in\Gamma(End\,TM)$ defines the {\it Schouten
concomitant} \cite{MM}
\begin{equation}\label{SchoutenR}R(\pi,A)(\alpha,X) =
\sharp_\pi(L_X(\alpha\circ A))- (L_{\sharp_\pi\alpha}A)(X)-
\sharp_\pi(L_{AX}\alpha)\end{equation} (note that our sign
convention for $\sharp_\pi$ and $\flat_\sigma$ is opposite to that
of \cite{MM}), which is equivalent to the $T^*M$-valued bivector
field \cite{VN}
\begin{equation}\label{Schouten} C_{(\pi,A)}(\alpha,\beta) =
\beta\circ L_{\sharp_\pi\alpha}A - \alpha\circ L_{\sharp_\pi\beta}A
+d(\pi(\alpha,\beta))\circ A - d(\pi(\alpha\circ A,\beta))
\end{equation} in the sense that
$$<R(\pi,A)(\alpha,X),\beta> =-
<C_{(\pi,A)}(\alpha,\beta),X>.$$ If the expression
(\ref{crosetforme}) of the bracket of $1$-forms is inserted in
(\ref{eqdinii}), it follows that condition ii) may be reformulated
as\\

\emph {ii') the Schouten concomitant $C_{(\pi,A)}$ vanishes}.\\

This interpretation of condition ii) has interesting consequences.
\begin{prop}\label{corespintegrabila} If the $2$-form $\sigma$ of
{\rm(\ref{matriceaPhi})} is symplectic, the structure $\Phi$ is
integrable iff the pair $(A,\sigma)$ is a symplectic-Nijenhuis
structure.\end{prop}
\begin{proof} We recall that a pair $(w\in\chi^2(M),A\in End\,TM)$
is a {\it Poisson-Nijenhuis structure} if $A$ and $w$ are
compatible, the Schouten-Nijenhuis bracket $[w,w]=0$, the Nijenhuis
tensor $\mathcal{N}_A=0$ and the Schouten concomitant $R(w,A)=0$
\cite{{MM},{VN}}. In particular, if $w$ comes from a symplectic form
$\sigma$, i.e., $\sharp_w\circ\flat_\sigma=-Id$, $(\sigma,A)$ is
called a {\it symplectic-Nijenhuis structure}. One can prove that
this happens iff the associated $2$-form $\sigma_A$ is also closed
\cite{{MM},{V2}}. The fundamental property of a Poisson-Nijenhuis structure
$(w,A)$ is the existence of a corresponding family of pairwise
compatible Poisson-Nijenhuis structures $(W,P_1(A))$, where
$\sharp_W=P_2(A)\circ\sharp_w$ and $P_{1,2}(A)$ are either
polynomials or convergent power series with constant coefficients
in the argument $A$, called the {\it Poisson hierarchy}
\cite{{MM},{VN}}. (The compatibility of two Poisson structures
$(w,W)$ means that the Schouten-Nijenhuis bracket $[w,W]=0$,
equivalently, that $w+cW$ $(c=const.)$ is again a Poisson bivector
field.)

Now, we notice that, for a non degenerate $2$-form $\sigma$,
conditions (\ref{2-reddezv}) are equivalent with
\begin{equation}\label{consec2dezv}
\sharp_\pi=(A^2-\epsilon Id)\circ\sharp_w,\;\sigma(AX,Y)=
\sigma(X,AY)\end{equation} (the condition for $\pi$ in (\ref{2-reddezv}) is
a consequence of (\ref{consec2dezv})) and we have
\begin{equation}\label{corespAPhi}\Phi=\left(\begin{array}{cc}
A&(A^2-\epsilon Id)\circ\sharp_w\vspace{2mm}\\
\flat_\sigma&-^t\hspace{-1pt}A\end{array}\right). \end{equation}
Then, if $d\sigma=0$ and $\Phi$ of (\ref{corespAPhi}) is integrable,
the integrability conditions ii'), iii), iv) show that $(\sigma,A)$
is a symplectic-Nijenhuis structure.

Conversely, if $(\sigma,A)$ is a symplectic-Nijenhuis structure
conditions ii'), iii), iv) for $\Phi$ hold and the integrability
condition i) follows from the Poisson hierarchy theorem.\end{proof}

Thus, the generalized, c.p.s. structures with a symplectic form
$\sigma$ are equivalent with the symplectic-Nijenhuis structures
with the same form $\sigma$. Moreover, a symplectic-Nijenhuis
manifold $(M,\sigma,A)$ is endowed with families of generalized,
c.p.s. structures defined by replacing $A$ by $P_1(A)$, where
$P_1(A)$ is either a polynomial or a convergent power series with
constant coefficients, in formula (\ref{corespAPhi}). Formula
(\ref{corespAPhi}) also shows that the Poisson structure $\pi$ of
$\Phi$ is compatible with the Poisson structure $w$ defined by the
symplectic form $\sigma$.

Furthermore, it is known that the symplectic-Nijenhuis structures
$(\sigma,A)$ of a manifold $M$ are in a bijective correspondence
with the compatible pairs of Poisson structures $(w,W)$
\cite{{MM},{V2}}. This correspondence sends $(\sigma,A)$ to $(w,W)$
where $\sharp_W=A\circ\sharp_w$, i.e., $W$ is the first new
Poisson structure of the Poisson hierarchy of $(\sigma,A)$.
Conversely, the pair $(w,W)$ is sent to
$(\sigma,-\sharp_W\circ\flat_\sigma)$.

This proves the following result
\begin{prop}\label{corolW} On a symplectic manifold $(M,\sigma)$
there exists a bijective correspondence between the Poisson
structures $W$ on $M$ that are compatible with $w$ and the
generalized, c.p.s. structures of $M$ which have the form $\sigma$
in their matrix representation. This correspondence is given by
\begin{equation}\label{corespcuW} W\mapsto
\Phi_W=\left(\begin{array}{cc} B&\sharp_\kappa\vspace{2mm}\\
\flat_\sigma&-^t\hspace{-1pt}B\end{array}\right)=
\left(\begin{array}{cc}
-\sharp_W\circ\flat_\sigma&-\epsilon\sharp_w-\sharp_W\circ\flat_\sigma
\circ\sharp_W\vspace{2mm}\\
\flat_\sigma&\flat_\sigma\circ\sharp_W\end{array}\right).
\end{equation} \end{prop}
\begin{rem} {\rm In the complex case,
the $\pm i$-eigenbundles of $\Phi_W$ are $L,\bar L$ where
\begin{equation}\label{definitiaL}
L=graph(\sharp_{-(W+iw)}:T^*_cM\rightarrow T_cM) \end{equation} and
$\bar L$ is the complex conjugate bundle of $L$. Indeed, the
conditions $[W,W]=0, [W,w]=0$ imply that $W+iw$ is a complex-valued
Poisson bivector field on $M$, hence $L$ defined by
(\ref{definitiaL}) is a complex Dirac structure. Furthermore, since
$\sigma$ is non degenerate, $L\cap\bar L=\{0\}$ and there exists a
unique, generalized, complex structure with the $\pm i$-eigenbundles
$L,\bar L$. In order to show that this structure is $\Phi_W$ it
suffice to compute the $i$-eigencomponent of $(X,\alpha)\in
T^{big}M$ with respect to $\Phi_W$:
\begin{equation}\label{icomponenta}
\frac{1}{2}(Id-i\Phi_W)(X,\alpha)=\frac{1}{2}
(X-i(BX+\sharp_\kappa\alpha),\alpha-i(\flat_\sigma X-\alpha\circ B))
\end{equation} $$=\frac{1}{2}(\sharp_{-(W+iw)}(\alpha -i(\flat_\sigma X
-\alpha\circ B)),\alpha -i(\flat_\sigma X -\alpha\circ B)).$$
Similar computations of eigenbundles may be done in the paracomplex
and subtangent cases.}\end{rem}

Other connections with Poisson-Nijenhuis structures are given by
\begin{prop}\label{corolPsNij} a) If $\Phi$ is integrable and
$\sigma$ is closed, $(\pi,A)$ is a Poisson-Nijenhuis structure on
$M$. b) If $(\pi,A)$ is a symplectic-Nijenhuis structure, $\Phi$ is
integrable iff the forms $\sigma$ and $\sigma_A$ are closed. c) Let
$\Phi$ be a generalized, almost c.p.s. structure on $M$ such that
$\pi$ is a Poisson bivector field, $A$ is a Nijenhuis tensor field,
and the $2$-forms $\sigma,\sigma_A$ are closed. Then $\Phi$ is
integrable.
\end{prop}
\begin{proof} Assertions a) and b) are trivial. For c), the
only condition we still have to check is $C_{(\pi,A)}=0$. For this
purpose, we write down the following formula, which holds for any
tensor fields $\pi\in\chi^2(M),\sigma\in\Omega^2(M)$ and is
equivalent with formula (B.3.9) of \cite{MM},
\begin{equation}\label{Magri1} C_{(\pi,\sharp_\pi\circ\flat_\sigma)}
(\alpha,\beta) = i(\sharp_\pi\beta)i(\sharp_\pi\alpha)d\sigma
-(i(\beta)i(\alpha)[\pi,\pi])\circ\flat_\sigma,
\end{equation} where $[\pi,\pi]$ is the Schouten-Nijenhuis
bracket. For a generalized, almost c.p.s. structure $\Phi$, the
first condition (\ref{2-reddezv}) changes (\ref{Magri1}) to
\begin{equation}\label{Magri1'} -C_{(\pi,A^2)} =
i(\sharp_\pi\beta)i(\sharp_\pi\alpha)d\sigma-
(i(\beta)i(\alpha)[\pi,\pi])\circ\flat_\sigma. \end{equation} Under
the hypotheses of iii), (\ref{Magri1'}) gives $C_{(\pi,A^2)}=0$, and
$(\pi, A^2)$ is a Poisson-Nijenhuis structure. Then, by the {\it
hierarchy theorem} for Poisson-Nijenhuis structures
\cite{{MM},{VN}}, $(\pi,A)$ is also a Poisson-Nijenhuis structure
and $C_{(\pi,A)}=0$.\end{proof}
\begin{rem}\label{obstwisted} {\rm In \cite{LMTZ} it
is shown that supersymmetry is also related with generalized, almost
complex structures that are integrable (i.e., have a vanishing
tensor (\ref{tensNij})) with respect to the \v{S}evera-Weinstein
Courant bracket \cite{SW}
\begin{equation}\label{crosetSWC}
[(X,\alpha),(Y,\beta)]=([X,Y],L_X\beta-L_Y\alpha\end{equation}
$$+\frac{1}{2}d(\alpha(Y)-\beta(X)) -i(Y)i(X)\Lambda),$$ where
$\Lambda$ is a closed $3$-form on $M$, and these new integrability
conditions are expressed in local coordinates. The computations done
to prove Theorem \ref{thCrainic} may be easily extended to
(\ref{crosetSWC}), and it follows that $\mathcal{N}_\Phi$ with
brackets (\ref{crosetSWC}) is zero iff

i) the bivector field $\pi$ defines a Poisson structure on $M$;

ii) the Schouten concomitant of the pair $(\pi,A)$ satisfies the
condition \begin{equation}\label{condScptSWC}
C_{(\pi,A)}(\alpha,\beta)=i(\sharp_\pi\beta)i(\sharp_\pi\alpha)
\Lambda;\end{equation}

iii) the Nijenhuis tensor of $A$ is given by the formula
\begin{equation}\label{Nijptintegrab'} \mathcal{N}_A(X,Y) =
\sharp_\pi[i(Y)i(X)d\sigma + i(AY)i(X)\Lambda-i(AX)i(Y)\Lambda];
\end{equation}

iv) the exterior differential of the {\it associated form}
$\sigma_A$ satisfies the equality
\begin{equation}\label{difsigmaA'}
d\sigma_A(X,Y,Z)-\epsilon\Lambda(X,Y,Z)\end{equation}} $$=
\sum_{Cycl(X,Y,Z)}[d\sigma(AX,Y,Z)+\Lambda(AX,AY,Z)].$$
\end{rem}

Following Crainic \cite{Cr}, we will say that a generalized, almost
c.p.s. structure $\Phi$ is {\it non degenerate} if its bivector
field $\pi$ is non degenerate. Then, we will denote by $\varpi$ the
non degenerate $2$-form defined by
$\flat_\varpi\circ\sharp_\pi=-Id,$ and (\ref{2-reddezv}) implies
\begin{equation}\label{sigmaptHitchin}
\flat_\sigma=\flat_\varpi\circ A^2-\epsilon\flat_\varpi.
\end{equation}
\begin{theorem}\label{integrinndeg} {\rm\cite{Cr}}
Let $\Phi$ be a non degenerate, generalized, almost c.p.s.
structure. Then, $\Phi$ is integrable iff $\pi$ is Poisson and the
$2$-form $\varpi_A$ is closed.\end{theorem}

A pair $(\varpi,A)$ where $\varpi$ is a symplectic form and $A$ is a
compatible $(1,1)$-tensor field is called a {\it Hitchin pair} if
the associated $2$-form $\varpi_A$ is closed \cite{Cr}. The previous
theorem may be used to show the existence of a $1-1$ correspondence
between each of the three classes of non degenerate, integrable,
almost c.p.s. structures (separately) and Hitchin pairs, which is
defined by the formula \cite{Cr},
\begin{equation}\label{corespondenta}
(\varpi,A)\mapsto \left(\begin{array}{cc}A&\sharp_\pi\vspace{2mm}\\
\flat_{\varpi_{A^2}}-\epsilon\flat_\varpi&-\,^t\hspace{-1pt}A\end{array}\right).
\end{equation}
Accordingly, Crainic's results on Lie groupoids and algebroids
connected with generalized complex structures have corresponding
variants for generalized paracomplex and tangent structures.\\

We continue the presentation of the basic results on generalized
c.p.s. structures by indicating some examples.

\begin{example}\label{exemplul2} {\rm \cite{Galt} For any classical c.p.s.
structure $A$ on $TM$, the matrix
\begin{equation}\label{eqex2}  \left(\begin{array}{cc}
A&0\vspace{2mm}\\0&-\,^t\hspace{-1pt}A\end{array}\right)
\end{equation}
is a generalized c.p.s. structure, respectively.}\end{example}
\begin{example}\label{exemplul1} {\rm \cite{{Galt},{Cr}}
Any symplectic form $\varpi$ produces the c.p.s. structures
\begin{equation}\label{eqex1} \left(\begin{array}{cc}0&\sharp_\pi\vspace{2mm}\\
-\epsilon\flat_\varpi&0\end{array}\right),\hspace{5mm}
\left(\begin{array}{cc}Id&\sharp_\pi\vspace{2mm}\\
(1-\epsilon)\flat_\varpi&-Id\end{array}\right), \end{equation} where
$\flat_\varpi\circ\sharp_\pi=-Id$, associated with the Hitchin pairs
$(\varpi,0),(\varpi,Id)$, respectively.}\end{example}
\begin{example}\label{exemplul4} {\rm \cite{BB} If $(M,F)$ is a locally
product manifold with {\it structural foliations} $ \mathcal{F}_1,
\mathcal{F}_2$ (i.e., $TM=T\mathcal{F}_1\oplus T\mathcal{F}_2$), and
if these foliations have generalized c.p.s. structures
$\Phi_1,\Phi_2$ along the leaves, which are differentiable on $M$,
then $\Phi=\Phi_1\oplus \Phi_2$ is a generalized c.p.s. structure on
$M$. Sometimes, it is interesting to change $\Phi_1$ by its {\it
opposite structure} $^{op}\hspace{-1pt}\Phi_1$, which is defined by
changing the sign of the tensor fields $\pi,\sigma$ in the matrix
(\ref{matriceaPhi}) of $\Phi_1$, and use the {\it twisted direct
sum} $^{op}\hspace{-1pt}\Phi_1\oplus \Phi_2$. Theorem
\ref{thCrainic} shows that $\Phi$ and $^{op}\hspace{-1pt}\Phi$ are
simultaneously integrable.}\end{example}
\begin{example}\label{exemplul5} {\rm Let $M$ be a complex
analytic manifold. Then we may define the notion of a {\it
holomorphic Dirac structure} in the same way as a real Dirac
structure, using the holomorphic tangent and cotangent bundles of
$M$. If $L_1,L_2$ are two holomorphic Dirac structures on $M$,
$L_1\oplus\bar L_2$ (where the bar denotes complex conjugation)
obviously is the $\sqrt{-1}$-eigenbundle of a generalized, complex
structure of $M$. (Hitchin's example of the generalized, complex
structure associated with a holomorphic Poisson structure of $M$
\cite{Ht2} is a particular case of the previous construction.)}
\end{example}

Another basic notion of the theory is that of gauge equivalence.
This notion is based on Hitchin's remark \cite{Ht1} that, for any
closed $2$-form $B$ on $M$, the mapping
\begin{equation}\label{gaugetr} (X,\alpha) \mapsto
(X,\alpha+i(X)B),\end{equation} called a $B$-{\it field} or {\it
gauge} transformation (equivalence) is a bundle automorphism $
\mathcal{B}$ of $T^{big}M$ which is compatible with the metric $g$
and the Courant bracket (\ref{crosetC}). The matrix of the gauge
transformation (\ref{gaugetr}) is
\begin{equation}\label{gaugetr2}
\mathcal{B}=\left( \begin{array}{cc}Id&0\vspace{2mm}\\
\flat_B&Id\end{array}\right),\end{equation} and $ \mathcal{B}$ acts
on generalized, almost c.p.s. structures by the invertible mapping
$\Phi\mapsto \mathcal{B}^{-1}\Phi\mathcal{B}$. By Proposition
\ref{integrprinDirac}, in the complex and paracomplex cases $
\mathcal{B}$ also preserves integrability.

From the algebraic point of view, $B$-field transformations may be
defined in the same way for $E\oplus E^*$, where $E$ is an arbitrary
vector bundle with a generalized, almost c.p.s. structure and
$B\in\Gamma\wedge^2E^*$ but, of course, no properties of any bracket
will be involved.

The action of $\mathcal{B}$ sends the matrix (\ref{matriceaPhi}) of
$\Phi$ to the matrix
\begin{equation}\label{gaugepePhi}
\left(\begin{array}{cc}
A+\sharp_\pi\circ\flat_B&\sharp_\pi\vspace{2mm}\\
\flat_\sigma - \flat_B\circ\sharp_\pi\circ\flat_B -\flat_B\circ A
-\,^t\hspace{-1pt}A\circ\flat_B&-\,^t\hspace{-1pt}A -
\flat_B\circ\sharp_\pi\end{array}\right).\end{equation} Thus, the
Poisson bivector field of $\Phi$ is preserved, the tensor field $A$
goes to $A+\sharp_\pi\circ\flat_B$, and the $2$-form $\sigma$ is
changed to
\begin{equation}\label{gaugesigma} \sigma'(X,Y)=\sigma(X,Y) +
\pi(\flat_BX,\flat_BY) - B(AX,Y) - B(X,AY).\end{equation}
\begin{example}\label{exdegauge} {\rm Any non degenerate c.p.s.
structure $\Phi$ represented by the right hand side of the mapping
(\ref{corespondenta}) is gauge equivalent with the symplectic
structure $\varpi$ seen as the first matrix (\ref{eqex1}), by the
field $B=\varpi_A$, and seen as the second matrix (\ref{eqex1}), by
the field $B=\varpi_A-\varpi$.} \end{example}
\begin{example}\label{ex2gauge} {\rm The
$B$-transform of the structure (\ref{eqex2}) is of the form
\begin{equation}\label{eqex3}  \left(\begin{array}{cc}
A&0\vspace{2mm}\\
\flat_\sigma&-\,^t\hspace{-1pt}A\end{array}\right),
\end{equation}
where $$\sigma(X,Y)=-B(AX,Y)-B(X,AY).$$ Conversely, if the Poisson
structure of a generalized c.p.s. structure $\Phi$ is zero, the
structure is gauge equivalent with a classical c.p.s. structure,
seen as (\ref{eqex2}), iff there exists a closed $2$-form $B$ such
that \begin{equation}\label{Bcuclasic}
\sigma(X,Y)=B(AX,Y)+B(X,AY).\end{equation} For $\epsilon=\pm1$, the
general algebraic solution of (\ref{Bcuclasic}) is
\begin{equation}\label{Bptclasic} B(X,Y) =
\frac{\epsilon}{2}\sigma(AX,Y)+B'(X,Y), \end{equation} where
$B'(AX,Y)+B'(X,AY)=0$, and we see that $\Phi$ defined by
(\ref{eqex3}), where $\sigma_A$ is closed, is gauge equivalent with
a classical structure.}
\end{example}

The importance of gauge equivalence is shown by the local structure
theorems of Gualtieri \cite{Galt} and Abouzaid-Boyarchenko
\cite{AB}, which show that any generalized, complex manifold is
gauge equivalent with the direct sum of a symplectic and a classical
complex structure in a neighborhood of a point. One also has the
algebraic result that any generalized, complex structure of a vector
space is gauge equivalent with such a direct sum \cite{BB}. We
extend this algebraic result in the following proposition.
\begin{prop} \label{propBB} Let $(M,\Phi)$ be a generalized,
almost c.p.s. manifold, where $\Phi$ is defined by the matrix
{\rm(\ref{matriceaPhi})} and the bivector field $\pi$ is Poisson.
Let $S$ be a symplectic leaf  of $\pi$. Let $\nu S$ be a normal
bundle of $S$, i.e., a subbundle of $T_SM$ such that
\begin{equation}\label{descloc}
T_SM=TS\oplus \nu S.\end{equation} Then, the restriction of $\Phi$
to the bundle $T^{big}_SM$ is algebraically gauge equivalent with a
direct sum of a symplectic structure on $T^{big}S$ and a c.p.s.
structure on $\nu^{big}S=\nu S\oplus\nu^*S$.
\end{prop}
\begin{proof} In the conclusion, by a symplectic structure we mean
a structure defined by the first matrix (\ref{eqex1}) and by a
c.p.s. structure we mean a matrix (\ref{eqex2}) which satisfies the
algebraic conditions of a generalized c.p.s. structure.

First, we show the existence of $B\in\Gamma\wedge^2T^*_SM$ such that
the $B$-field equivalent structure
$$\Phi'= \left(\begin{array}{cc} A'&\sharp_\pi\vspace{2mm}\\
\flat_{\sigma'}&-\,^t\hspace{-1pt}A'\end{array}\right)$$ of $\Phi$
has the property that $\nu S$ invariant by $A'$. Indeed, by
(\ref{gaugepePhi}), this means that we have to choose $B$ such that
\begin{equation}\label{ptBauxiliar}
\sharp_\pi\flat_B(V)=-pr_{TS}A(V),\hspace{3mm}\forall V\in \nu
S,\end{equation} where the projection is defined by the
decomposition (\ref{descloc}). Since $S$ is a symplectic leaf of
$\pi$, there exists a unique $\lambda\in T^*S$ such that
$pr_{TS}\circ A(V)=\sharp_\pi\lambda$, and $V\mapsto-\lambda$ yields
a well defined mapping $\varphi:\nu S\rightarrow T^*S$ such that
(\ref{ptBauxiliar}) is satisfied if $\flat_B|_{\nu S}=\varphi$. This
mapping extends to a mapping $\flat_B: TS\oplus \nu S\rightarrow
T^*S\oplus \nu^*S^*$ defined in matrix form by
$$\flat_B= \left(\begin{array}{cc} 0&\varphi\vspace{2mm}\\
-\,^t\hspace{-1pt}\varphi&0\end{array}\right),$$ which is associated
with a $2$-form $B$. (Above and hereafter $T^*S,\nu^*S$ are seen as
the terms of the decomposition $T^*_SM = T^*S\oplus\nu^*S$ induced
by (\ref{descloc}).)

Furthermore, we shall see that $TS\oplus T^*S$ and $\nu S\oplus
\nu^*S$ are invariant by $\Phi'$ and the latter is the direct sum of
its restrictions to these invariant subbundles, which are of the
form
\begin{equation}\label{Phirestrans}
\Phi'|_{TS\oplus T^*S} = \left( \begin{array}{cc}A'&\sharp_\pi\vspace{2mm}\\
\flat_{\sigma'}&-\,^t\hspace{-1pt}A'\end{array}\right),\;
\Phi'|_{\nu S\oplus \nu^*S} = \left( \begin{array}{cc}A'&0\vspace{2mm}\\
\flat_{\sigma'}&-\,^t\hspace{-1pt}A'\end{array}\right).
\end{equation} Indeed, we obviously have $A'=A'|_{TS}\oplus A'|_{\nu
S}$ and $\sharp_\pi=(\sharp_\pi)|_{T^*S}\oplus0$ (remember that
$TS=im\,\sharp_\pi$ and $\nu^*S=ann\,TS$). In order to see that
$\flat_{\sigma'} =
(\flat_{\sigma'})|_{TS}\oplus(\flat_{\sigma'})|_{\nu S}$ we have to
check the $\sigma'$-orthogonality of $TS$ and $\nu S$, which is seen
as follows. For $X=\sharp_\pi\xi\in TS$ $(\xi\in T^*S)$ and $V\in\nu
S$, (\ref{2-reddezv}) for $\Phi'$ implies
$$\sigma'(V,X)=<\flat_{\sigma'}V,\sharp_\pi\xi> =
-<\sharp_\pi\flat_{\sigma'}V,\xi> $$ $$=-<\epsilon V-A'^2V,\xi>=0.$$

Now, since $\Phi'|_{TS\oplus T^*S}$ is non degenerate, Example
\ref{exdegauge} tells us that, algebraically, this component of
$\Phi'$ is gauge equivalent with the symplectic structure defined by
the first matrix (\ref{eqex1}) associated to $\sharp_\pi|_{T^*S}$.

Then, in view of (\ref{Bcuclasic}), a gauge transformation that
sends $\Phi'|_{\nu S\oplus \nu^*S}$ to a structure with a matrix
form (\ref{eqex2}) is defined by a new $B\in\Gamma\wedge^2\nu^*S$
such that
$$\sigma'(V_1,V_2)=B(A'V_1,V_2)+B(V_1,A'V_2)$$ $(V_1,V_2\in\nu S)$.
This condition holds if we ask
\begin{equation}\label{Baux2}
B(A'V_1,V_2)=\frac{1}{2}\sigma'(V_1,V_2).\end{equation} Such a form
$B$ exists: for $\epsilon=\pm1$ $A'|_{\nu S}$ is non degenerate, and
(\ref{Baux2}) fully defines $B$; for $\epsilon=0$, (\ref{Baux2})
defines $B$ on $im\,A'|_{\nu S})$ and we may use an arbitrary
extension to $\nu S$.

The composition of all the algebraic gauge transformations described
above yields the required conclusion. \end{proof}
\begin{rem}\label{obsbetafield} {\rm From the algebraic point of
view again, it is also interesting to refer to a dual notion of
$\beta$-{\it field transformation} \cite{{Galt},{BB}}
\begin{equation}\label{gaugetransf}
\Phi\mapsto\mathcal{B}^{'-1}\Phi\mathcal{B}',\end{equation} where
$$ \mathcal{B}'=\left(\begin{array}{cc}
Id&\sharp_\beta\vspace{2mm}\\
0&Id\end{array}\right)\hspace{5mm}(\beta\in\chi^2(M)).$$
Preservation of integrability by a $\beta$-field transformation is
rare. An example is given by formula (\ref{corespcuW}) where
$\Phi_W$ is the result of the $W$-field transformation of the
symplectic structure $\sigma$ seen as in the first matrix
(\ref{eqex1}) and $\Phi_W$ is integrable.}
\end{rem}

We finish this section by referring to a notion of {\it generalized
c.p.s. mapping} (generalized, holomorphic mapping in the complex
case). This is not simple because both contravariant and covariant
tensor fields are involved. The most appropriate definition seems to
be that of Crainic \cite{Cr}, even though it is very restrictive. We
justify Crainic's definition as follows.

A mapping $f:(M_1,\Phi_1)\rightarrow (M_2,\Phi_2)$, where
$\Phi_1,\Phi_2$ are generalized c.p.s. structures, produces
relations
\begin{equation}\label{frelatie} f^{rel}_x =
\{((X,f^*\alpha),(f_*X,\alpha))\,/\,X\in T_xM_1,\alpha\in
T^*_xM_2\}\end{equation} $$\subseteq T_x^{big}M_1\times
T_{f(x)}^{big}M_2,$$ defined $\forall x\in M_1$. The mapping $f$
will be called a generalized c.p.s. mapping if, $\forall x\in M_1$,
$\forall((X,f^*\alpha), (f_*X,\alpha))\in f_x^{rel}$ one has
$(\Phi_1(X,f^*\alpha),\Phi_2(f_*X,\alpha))\in f_x^{rel}$.

Accordingly, if the matrices of $\Phi_1,\Phi_2$ are as in
(\ref{matriceaPhi}) with indices $1,2$, respectively, $f$ is
generalized c.p.s. iff, $\forall x\in M_1$ and $\forall X\in
T_xM_1,\forall\alpha\in T^*_{f(x)}M_2$ one has
\begin{equation}\label{fgenhol} \begin{array}{lcl}
A_2(f_*X)+\sharp_{\pi_2}\alpha&=&f_*(A_1X+\sharp_{\pi_1}f^*\alpha),
\vspace{2mm}\\ \flat_{\sigma_1}X-(f^*\alpha)\circ A_1& = &
f^*(\flat_{\sigma_2}(f_*X) - \alpha\circ A_2).\end{array}
\end{equation}

Furthermore, if we look at (\ref{fgenhol}) for either $X=0$ or
$\alpha=0$, we see that $f$ is generalized c.p.s. iff the following
three conditions required in \cite{Cr} hold
\begin{equation}\label{fgenhol1}
\pi_2=f_*\pi_1,\;\sigma_1=f^*\sigma_2,\;A_2\circ f_*=f_*\circ A_1.
\end{equation}

Following is an example of utilization of the notion of a c.p.s.
mapping which shows its restrictive character.

We define a {\it generalized c.p.s. Lie group} to be a real Lie
group $G$ endowed with a generalized c.p.s. structure $\Phi$ such
that the multiplication mapping $\mu(g_1,g_2)=g_1g_2$ is a
generalized c.p.s. mapping $\mu:(G\times G,\Phi\oplus\Phi)
\rightarrow (G,\Phi)$.

If the matrix of $\Phi$ is that of formula (\ref{matriceaPhi}), the
first condition (\ref{fgenhol1}) tells that the pair $(G,\pi)$ is a
Poisson-Lie group \cite{V-carte}. The second condition
(\ref{fgenhol1}) gives
\begin{equation}\label{conddoua} \sigma_{g_1}(X,X') + \sigma
_{g_2}(Y,Y') = \sigma_{g_1g_2}(L_{g_1^*}Y + R_{g_2*}X, L_{g_1^*}Y' +
R_{g_2*}X'),\end{equation} $\forall g_1,g_2\in G,\, X,X'\in
T_{g_1}G,\, Y,Y'\in T_{g_2}G$ and where $L,R$ denote left and right
translations, respectively. Indeed, it is easy to see (e.g.,
\cite{V-carte}) that

$$\mu_*(X,Y)=L_{g_1*}Y+ R_{g_2*}X.$$ For $X=X'=0$, respectively, $Y=Y'=0$
condition (\ref{conddoua}) shows that $\sigma$ is left-invariant and
right-invariant, respectively. Then, the case $X'= Y=0$ shows that
the only possibility is $\sigma=0$. Finally, in a similar way, we
see that the third condition (\ref{fgenhol1}) simply means that $A$
is a bi-invariant tensor field on $G$. Moreover, in view of
$\sigma=0$ and of (\ref{2-reddezv}), in the complex and paracomplex
cases $A$ is an almost c.p. structure, respectively.

Now, if we also look at the integrability conditions given by
Theorem \ref{thCrainic}, we get
\begin{prop} \label{propgrpgenc} A generalized complex, respectively
paracomplex, Lie group is a classical complex, respectively
paracomplex Lie group $(G,A)$ endowed with a multiplicative Poisson
bivector field $\pi$ such that the pair $(\pi,A)$ is a
Poisson-Nijenhuis structure.
\end{prop}
\section{Reduction of generalized c.p.s. structures}
Reduction theory is an important chapter of symplectic and Poisson
geometry. {\it Geometric reduction} leads to a symplectic,
respectively Poisson, structure on a quotient of a submanifold of a
given symplectic or Poisson manifold. If the submanifold is obtained
as a non critical level set of a momentum map of a Hamiltonian group
action, one has the {\it Marsden-Weinstein reduction}, which has
many applications in mechanics and physics.

The general, Poisson framework of geometric reduction was given by
Marsden and Ratiu \cite{MRt} and we briefly describe it as follows.
Let $(M,\pi)$ be a Poisson manifold, and $\iota:N\hookrightarrow M$
a submanifold. A subbundle $E\subseteq T_NM$ is a {\it
reduction-control bundle} on $N$ if a) $E\cap TN=T\mathcal{F}$ where
$ \mathcal{F}$ is a foliation of $N$ by the fibers of a submersion
$s:N\rightarrow Q$, b) $\forall\varphi,\psi\in C^\infty(M)$ such
that $d\varphi|_N,d\psi|_N\in ann\,E$ the Poisson bracket
$\{\varphi,\psi\}$ satisfies the same condition
$d\{\varphi,\psi\}_\pi|_N\in ann\,E$ ($ann$ denotes the annihilator
of the bundle $E$), c) $\sharp_\pi(ann\,E) \subseteq TN+E$. The
Marsden-Ratiu reduction theorem says that if $E$ is a
reduction-control bundle there exists a unique Poisson structure
$\pi_Q$ on $Q$ such that $\forall\lambda\in T^*Q$ one has
\begin{equation}\label{piredus} \sharp_{\pi_Q}\lambda=
s_*(pr_{TN}\sharp_\pi \widetilde{s^*\lambda}),\end{equation} where
$\widetilde{s^*\lambda}$ is an extension of $s^*\lambda$ to $T_NM$
such that $\widetilde{s^*\lambda}|_N\in ann\,E$ (an $E$-{\it
controlled extension}). Formula (\ref{piredus}) holds $\forall x\in
N,y=s(x)\in Q$. The projection $pr_{TN}$ is defined by the
decomposition of condition c) for $E$; it may not be uniquely
defined but, any two values of this projection differ by a vector in
$T\mathcal{F}$ and the projection by $s_*$ is well defined.
$E$-controlled extensions of $s^*\lambda$ may be obtained by asking
them to vanish on a normal bundle $\nu N$ of $N$ in $M$ (i.e.,
$T_NM=TN\oplus \nu N$) which is of the form $\nu N=E'\oplus C$,
where $E'$ is a complement of $T\mathcal{F}$ in $E$ and $C$ is a
complement of $E$ in $T_NM$. The independence of $\pi_Q$ on the
choice of the controlled extensions is part of the proof of the
reduction theorem \cite{MRt}. The Poisson structure $\pi_Q$ is the
{\it reduction of $\pi$ via $(N,E)$}.

Now, let us consider a generalized c.p.s. manifold $(M,\Phi)$, where
$\Phi$ has the matrix (\ref{matriceaPhi}), and a submanifold
$\iota:N\hookrightarrow M$ with a $\pi$-reduction-control bundle $E$
and with the reduced structure $\pi_Q$ of the quotient manifold
$Q=M/\mathcal{F}$ $(T\mathcal{F}=E\cap TN)$. We would like to be
able to reduce the whole structure $\Phi$ to $Q$, and we will prove
a theorem which shows that, if hypotheses that ensure the
reducibility of $\Phi$ to a generalized, almost c.p.s. structure are
satisfied, the reduced structure is integrable.
\begin{theorem}\label{thredgeom} Assume that the configuration
$(M,\Phi,N,E,Q,\pi_Q)$ considered above satisfies the following
hypotheses: 1) $\sharp_\pi(ann\,E)\subseteq TN$; 2) $ A(TN)\subseteq
TN,A(E)\subseteq E$ and $A|_{TN}$ sends $ \mathcal{F}$-projectable
vector fields $X\in\Gamma TN$ to $ \mathcal{F}$-projectable vector
fields $AX\in\Gamma TN$; 3) $\forall Z\in\Gamma E$ one has
$\iota^*(i(Z)\sigma)=0,\iota^*(i(Z)d\sigma)=0$. Then $Q$ has unique
tensor fields $A_Q,\sigma_Q$ which, together with the reduced
Poisson structure $\pi_Q$, define a generalized c.p.s. structure
$\Phi_Q$ on $Q$.\end{theorem}
\begin{proof} We know that $\pi_Q$ exists and is Poisson from the
Marsden-Ratiu theorem; hypothesis 1) means that $ann\,E$ and
$ann\,TN$ are $\pi$-orthogonal, it is stronger than property c) of
$E$ and is required for the continuation of the present proof. The
existence of a projection $A_Q$ of $A$ is obvious from hypothesis
2). Hypothesis 3) ensures that $\iota^*\sigma=s^*\sigma_Q$ for a
well defined $2$-form $\sigma_Q$; indeed, in foliation theory it is
known that the existence of $\sigma_Q$ is ensured by the conditions
$i(Z)(\iota^*\sigma)=0$, $L_Z(\iota^*\sigma)=0$, $\forall Z\in\Gamma
T\mathcal{F}$ \cite{Mol} and these conditions are implied by
hypothesis 3). Thus, the theorem will be proven if we check the
algebraic and integrability conditions of a generalized c.p.s.
structure, which is done as follows.
$$1.\hspace{1cm} A_Q^2[X]_{T_x\mathcal{F}}= [A^2X]_{T_x\mathcal{F}} =
[\epsilon X - \sharp_\pi\circ\flat_\sigma X]_{T_x\mathcal{F}} =
\epsilon [X]_{T_x\mathcal{F}} -
\sharp_{\pi_Q}\circ\flat_{\sigma_Q}[X]_{T_x\mathcal{F}},$$ where
$x\in N$, we have identified a vector in $T_{s(x)}Q$ with an
equivalence class $[X]_{T_x\mathcal{F}}$ modulo $T_x\mathcal{F}$
$X\in T_xN$ on $N$ and
$$\sharp_{\pi_Q}\circ\flat_{\sigma_Q}[X]_{T_x\mathcal{F}} =
[\sharp_\pi\circ\flat_\sigma X]_{T_x\mathcal{F}}$$ in view of
(\ref{piredus}) and because, by hypotheses 3), $\flat_\sigma\tilde
X$ is a controlled extension of
$s^*(\flat_{\sigma_Q}[X]_{T_x\mathcal{F}})$.\vspace{2mm}\\
2.\hspace{1cm}The compatibility of $\sigma_Q$ with $A_Q$ is trivial
and the compatibility of $\pi_Q$ with $A_Q$ is a consequence of
(\ref{piredus}) and of the fact that, if $\lambda\in\Omega^1(Q)$ and
$\tilde\lambda$ is a controlled extension of $s^*\lambda$,
$\tilde\lambda\circ A$ is a controlled
extension of $s^*(\lambda\circ A_Q)$.\vspace{2mm}\\
3.\hspace{1cm}The integrability condition i), of Theorem
\ref{thCrainic} holds, and checking condition iv) of Theorem
\ref{thCrainic} is trivial because we have
$s^*(\sigma_{Q_{A_Q}})=\iota^*\sigma_A$, and $s^*$ is
injective.\vspace{2mm}\\ 4.\hspace{1cm}Condition iii) of Theorem
\ref{thCrainic} holds since we have $$
\mathcal{N}_{A_Q}([X]_{T_x\mathcal{F}},[Y]_{T_x\mathcal{F}}) = [
\mathcal{N}_A(X,Y)]_{T_x\mathcal{F}}$$ $$=
[\sharp_\pi(i(Y)i(X)d\sigma)]_{T_x\mathcal{F}} =
\sharp_{\pi_Q}(i([Y]_{T_x\mathcal{F}})i([X]_{T_x\mathcal{F}})d\sigma_Q);$$
the last equality holds because hypothesis 3) implies that
$i(Y)i(X)d\sigma$ is a controlled extension of
$s^*(i([Y]_{T_x\mathcal{F}})i([X]_{T_x\mathcal{F}})d\sigma_Q)$.\vspace{2mm}\\
5.\hspace{1cm}The proof of the fact that the Schouten concomitant
$C_{(\pi_Q,A_Q)}$ vanishes appears in the proof of the reduction
theorem for Poisson-Nijenhuis structures \cite{VN1}, page
92-93.\end{proof}
\begin{rem}\label{redusexplicit} {\rm From facts included in the
proof of Theorem \ref{thredgeom} we get the following explicit
formula for the reduced structure $\Phi_Q$
\begin{equation}\label{eqPhiQ} \Phi_Q([X]_{T_x\mathcal{F}},\lambda)
= (s_*pr_{TN}(pr_{TM}\Phi(X,\widetilde{s^*\lambda})),
(s^*)^{-1}\iota^*(pr_{T^*M}\Phi(X,\widetilde{s^*\lambda}))),\end{equation}
where $X\in T_xN,\lambda\in T^*_{s(x)}Q,x\in N$ and
$\widetilde{s^*\lambda}$ is a controlled extension of $s^*\lambda$.}
\end{rem}

For any submanifold $\iota:N\hookrightarrow (M,\Phi)$, where $\Phi$
is a generalized, almost c.p.s. structure, one has the following
differentiable field of planes along $N$
\begin{equation}\label{enlargedimage} \nu^{ei}N= pr_{TM}
(\Phi(TN\oplus ann\,TN)) \end{equation} $$=
\{AX+\sharp_\pi\alpha\,/\,X\in TN,\alpha\in ann\,TN\} =A(TN) +
\sharp_\pi(ann\,TN).$$ The field $\nu^{ei}N$, which may not have a
constant dimension, will be called the {\it enlarged image field} of
$N$.

We assume that $\Phi$ is integrable and that $\nu^{ei}N$ is a vector
bundle over $N$ (i.e., all its planes are of the same dimension).
Then, we shall discuss conditions ensuring that $\Phi$ can be
reduced via $(N,\nu^{ei}N)$.

From the last equality (\ref{enlargedimage}) it follows that
\begin{equation}\label{annpseudo} ann(\nu^{ei}N) = ann(A(TN))\cap
(ann\,TN)^{\perp_\pi}\end{equation} and, as a consequence of
(\ref{annpseudo}), we get $\sharp_\pi ann(\nu^{ei}N) \subseteq TN$,
which is condition 1) of Theorem \ref{thredgeom}. Moreover,
hypothesis $A(E)\subseteq E$ of condition 2) of Theorem
\ref{thredgeom} is implied by $A(TN)\subseteq TN$. Indeed, if we ask
the latter, we also have $^t\hspace{-1pt}A(ann\,TN) \subseteq
ann\,TN$ and (\ref{enlargedimage}) yields $A(\nu^{ei}N)\subseteq
\nu^{ei}N$.

Condition b) of reduction-control is equivalent with
\begin{equation}\label{bechivalent} (L_{\tilde Z}\pi)|_{ann\,E}=0,
\end{equation} where
$\tilde Z$ is an extension of $Z\in\Gamma E$ to $M$ (evaluate the
Lie derivative (\ref{bechivalent}) on the arguments
$d\varphi|_N,d\psi|_N$ of condition b)) and for $E=\nu^{ei}N$
(\ref{bechivalent}) is equivalent to
\begin{equation}\label{condbptnuei} (L_{A\tilde
X}\pi)|_{ann(\nu^{ei}N)}=0,\;
(L_{\sharp_\pi\tilde\alpha}\pi)|_{ann(\nu^{ei}N)}=0,\end{equation}
where $\tilde X$ extends $X\in\Gamma TN$ and $\tilde\alpha$ extends
$\alpha\in\Gamma(ann\,TN)$. The second condition (\ref{condbptnuei})
always holds. To see this, we use the characterization of Poisson
structures via the Schouten-Nijenhuis bracket, $[\pi,\pi]=0$, and
the fact that $-[\pi,.]$ is the Lichnerowicz coboundary $ \sigma$,
which is the contravariant, exterior differential on the Poisson
manifold $(M,\pi)$ \cite{V-carte}. Since $\forall\lambda\in
ann(\nu^{ei}N)$ we have $\sharp_\pi\lambda\in TN$, the usual
formulas of the Lie derivative and of the contravariant, exterior
differential (e.g., \cite{V-carte}, formula (4.8)) yield
$$(L_{\sharp_\pi\tilde\alpha}\pi)(\lambda_1,\lambda_2) = -
\sigma(\pi) = [\pi,\pi](\alpha,\lambda_1,\lambda_2)=0.$$

In what follows, {\it invariant} always means $A$-invariant, and we
assume that this condition holds. Then
\begin{equation}\label{intersectps} \nu^{ei}N\cap TN =
\{AX+\sharp_\pi\xi\,/\,X\in TN, \xi\in ann\,TN,\,\sharp_\pi\xi\in
TN\}\end{equation}
$$= A(TN) +
\sharp_\pi((ann\,TN)\cap(ann\,TN)^{\perp_\pi}),$$ and we shall
compute the Lie bracket of vector fields (\ref{intersectps}) with
differentiable $X$ and $\xi$ (if any). Using integrability
conditions i) (which implies $[\sharp_\pi\xi,
\sharp_\pi\eta]=\sharp_\pi\{\xi,\eta\}_\pi$), ii') (under the form
$R(\pi,A)=0$, where $R$ is defined by (\ref{SchoutenR})) and iii) of
Theorem \ref{thCrainic}, we get
\begin{equation}\label{intercroset}
[AX+\sharp_\pi\xi,AY + \sharp_\pi\eta]= A([AX,Y] +
[X,AY]\end{equation} $$- A[X,Y] + [\sharp_\pi\xi,Y] -
[\sharp_\pi\eta,X]) + \sharp_\pi(i(Y)i(X)d\sigma +
\{\xi,\eta\}_\pi$$ $$- L_{AY}\xi + L_Y(\xi\circ A) + L_{AY}\eta -
L_X(\eta\circ A)).$$

The first term of (\ref{intercroset}) is of the form required by
(\ref{intersectps}). Since the left hand side and the first term of
the right hand side of (\ref{intercroset}) are in $TN$, so is the
second term of the right hand side. Thus, the only condition
required in order to ensure that the bracket (\ref{intercroset})
belongs to $\nu^{ei}N\cap TN$ is \begin{equation}\label{condpt2.7}
i(Y)i(X)d\sigma + \{\xi,\eta\}_\pi- L_{AY}\xi + L_Y(\xi\circ
A)\end{equation}
$$ + L_{AY}\eta -
L_X(\eta\circ A)\in ann\,TN.$$ Using (\ref{crosetforme}) to evaluate
$\{\xi,\eta\}_\pi$ on $V\in TN$ we get $0$ (Lie derivatives are to
be computed using extensions of vector fields and forms from $N$ to
$M$ and the result does not depend on the choice of the extension),
hence, $\{\xi,\eta\}_\pi\in ann\,TN$. Similarly, the evaluation of
the last four terms of (\ref{condpt2.7}) on $V\in TN$ shows that
each of them belongs to $ann\,TN$. Therefore, (\ref{condpt2.7})
reduces to
\begin{equation}\label{dsigmaNzero} \iota^*d\sigma=0.
\end{equation}

Thus, if (\ref{dsigmaNzero}) holds and if $\nu^{ei}N\cap TN$ is a
distribution of planes on $N$ that has a constant dimension and
local generators $AX+\sharp_\pi\xi$ $(X\in
TN,\xi\in(ann\,TN)\cap(ann\,TN)^{\perp_\pi})$ where $X,\xi$ are
differentiable then $\nu^{ei}N\cap TN$ is a foliation $
\mathcal{F}_N$ of the submanifold $N$. By (\ref{intersectps}), the
existence of the required local generators is ensured if we ask
$(ann\,TN)\cap(ann\,TN)^{\perp_\pi})$ to have a constant dimension
or, equivalently, if we ask that $dim(TN+\sharp_\pi ann\,TN)=const.$

Moreover, we can also prove that $A|_N$ sends $
\mathcal{F}_N$-foliated vector fields to $\mathcal{F}_N$-foliated
vector fields. Let $X\in\chi^1(N)$ be $\mathcal{F}_N$-foliated and
take $Y\in\Gamma T\mathcal{F}_N$. We have to check that $[Y,AX]\in
T\mathcal{F}_N$. If $Y=AV$ with $V\in\chi^1(N)$, $A$-invariance,
condition (\ref{dsigmaNzero}), and the integrability condition iii)
of Theorem \ref{thCrainic} yield
$$ [AV,AX] = A[AV,X]+A[V,AX]-A^2[V,X] + \sharp_\pi(i(X)i(V)d\sigma)
\in T\mathcal{F}_N.$$ If $Y=\sharp_\pi\xi\in \Gamma TN$ where
$\xi\in ann\,TN$, the integrability condition ii') of Theorem
\ref{thCrainic} and the expression (\ref{SchoutenR}) of the Schouten
concomitant $R_{(\pi,A)}$ yield
$$[\sharp_\pi\xi,AX]=(L_{\sharp_\pi\xi}A)(X) + A[\sharp_\pi\xi,X]
= \sharp_\pi(L_X(\xi\circ A) - L_{AX}\xi) + A[\sharp_\pi\xi,X]\in
T\mathcal{F}_N.$$

Continuing to keep the $A$-invariance condition enforced, let us see
the meaning of hypothesis 3) of Theorem \ref{thredgeom}, where we
look at arguments $Z=AX$ $(X\in TN$) and $Z=\sharp_\pi\alpha$
$(\alpha\in ann\,TN)$. The conditions for $\sigma$ are
\begin{equation}\label{3cuAsialpha} \sigma(AX,Y)=0,\;
\sigma(\sharp_\pi\alpha,Y)=0,\hspace{1cm}Y\in TN,\alpha\in ann\,TN.
\end{equation} The second condition (\ref{3cuAsialpha}) is equivalent
with the invariance of $TN$ by $\sharp_\pi\circ\flat_\sigma$ and, in
view of (\ref{2-reddezv}), this is ensured by the $A$-invariance of
$N$.

The condition 3) for $d\sigma$ means
\begin{equation}\label{3cudAsialpha} d\sigma(AX,Y_1,Y_2)=0,\;
d\sigma(\sharp_\pi\alpha,Y_1,Y_2)=0,\end{equation} where
$\;X,Y_1,Y_2\in TN,\alpha\in ann\,TN$. The first condition
(\ref{3cudAsialpha}) is implied by (\ref{dsigmaNzero}) and the
second condition is a consequence of integrability condition iii),
Theorem \ref{thCrainic}.

Accordingly, we get the following reduction theorem.
\begin{theorem}\label{a2-areducere} Let $(M,\Phi)$ be a
generalized c.p.s. manifold and $\iota:N\rightarrow M$ an
$A$-invariant submanifold. Assume that the following hypotheses are
satisfied: 1) $dim\,\nu^{ei}N=const.$, $dim(\nu^{ei}N\cap
TN)=const.$ and $dim(TN+\sharp_\pi ann\,TN)=const.$, 2) the $2$-form
$\sigma$ satisfies the conditions $\iota^*\sigma_A=0$,
$\iota^*d\sigma=0$, 3) the foliation $ \mathcal{F}_N$, which exists
because of 1) and 2), consists of the fibers of a submersion
$s:N\rightarrow Q$, 4) the underlying Poisson structure satisfies
the first condition {\rm(\ref{condbptnuei})}. Then, $Q$ has a
reduced generalized c.p.s. structure $\Phi_Q$.\end{theorem}
\begin{proof} The hypotheses and the previous analysis show that
$N$ and $E=\nu^{ei}N$ satisfy all the hypotheses of Theorem
\ref{thredgeom}. \end{proof}
\begin{rem}\label{obsdupath2.2} {\rm
If the Poisson structure $\pi$ of the generalized c.p.s. structure
$\Phi$ is zero (e.g., $\Phi$ is a classical c.p.s. structure) and
$N$ is an invariant submanifold then $(A|_N,0,\iota^*\sigma)$ is a
generalized c.p.s. structure $\Phi_N$ of $N$ and reductions provided
by Theorem \ref{thredgeom} are just the projection of $\Phi_N$ to a
space of leaves.}\end{rem}

Another interesting field of planes along $N$ is the pseudo-normal
field of the submanifold $N$ with respect to the Poisson structure
$\pi$, $\nu_\pi N=\sharp_\pi(ann\,TN)$ (the name pseudo-normal was
introduced in \cite{Vtg}). It follows immediately that
\begin{equation}\label{annpiann} ann(\nu_\pi N)=\{\lambda\in
T^*_NM\,/\,\sharp_\pi\lambda\in TN\},\end{equation} therefore,
$\sharp_\pi(ann\,\nu_\pi N)\subseteq TN$. Furthermore, for two
vector fields in $\nu_\pi N\cap TN$ that are of the form
$\sharp_\pi\lambda_1,\sharp_\pi\lambda_2$ where
$\lambda_1,\lambda_2$ are differentiable and belong to
$(ann\,TN)\cap(ann\,\nu_\pi N)$, the bracket necessarily belongs to
$TN$ and it is given by \begin{equation}\label{crosetlambda12}
[\sharp_\pi\lambda_1,\sharp\pi\lambda_2]=\sharp_\pi\{\lambda_1,
\lambda_2\}_\pi.\end{equation} It is easy to check that
$\{\lambda_1, \lambda_2\}_\pi\in ann\,TN$, hence, the bracket
(\ref{crosetlambda12}) also belongs to $\nu_\pi N$. Thus, if the
field $\nu_\pi N\cap TN$ consists of planes of the same dimension
and is locally spanned by vector fields $\sharp_\pi\lambda$ with
differentiable $1$-forms $\lambda\in(ann\,TN)\cap(ann\,\nu_\pi N)$
then this field is an involutive distribution and we have a
foliation $\mathcal{C}(N)$ of $N$ such that $T\mathcal{C}=\nu_\pi
N\cap TN$ (see also \cite{V-carte}, p.104).

This situation leads to one more reduction theorem:
\begin{theorem}\label{a3-areducere} Let $N$ be an
$A$-invariant submanifold of a generalized c.p.s. manifold
$(M,\Phi)$. If $dim\,\nu_\pi N=const.$ and $dim(\nu_\pi N\cap
TN)=const.$, $N$ has a foliation $\mathcal{C}(N)$ with tangent
bundle $T\mathcal{C}= \nu_\pi N\cap TN$, and in case the leaves of
$\mathcal{C}(N)$ are the fibers of a submersion $s:N\rightarrow Q$,
$Q$ has a reduced generalized c.p.s. structure $\Phi_Q$ of $\Phi$
via $(N,\nu_\pi N)$.\end{theorem}
\begin{proof} The hypotheses of the corollary
imply conditions 1), 2) and 3) of Theorem \ref{thredgeom}. In
particular, $\nu_\pi N\cap TN$ is spanned by vector fields
$\sharp_\pi\lambda$ with differentiable $1$-forms
$\lambda\in(ann\,TN)\cap(ann\,\nu_\pi N)$ because the constancy of
the dimensions of $\nu_\pi N$, $\nu_\pi N\cap TN$ implies
$dim(TN+\nu_\pi N)=const.$, therefore $dim(ann\,TN)\cap(ann\,\nu_\pi
N)=const.$ The projectability of $A$, the existence of the reduced
Poisson structure $\pi_Q$, and the conditions of hypothesis 3) were
proven during the proof of Theorem \ref{a2-areducere} (the
difference between the present situation and the one in Theorem
\ref{a2-areducere} is that the vectors $\sharp_\pi\xi$ with $\xi\in
(ann\,TN)\cap(ann\,TN)^{\perp_\pi}$ suffice to span
$T\mathcal{C}(N)$.) \end{proof}
\begin{corol}\label{a4-areducere} Let $(M,\Phi)$ be a non
degenerate, generalized c.p.s. manifold where $\Phi$ is associated
to the Hitchin pair $(\varpi,A)$. Let
$\iota:N\hookrightarrow(M,\Phi)$ be an $A$-invariant submanifold
such that: 1) $rank\,\iota^*\varpi=const.$, 2) the leaves of the
foliation $\mathcal{C}(N)$ $(T\mathcal{C}=\nu_\pi N\cap TN)$  are
the fibers of a submersion $s:N\rightarrow Q$. Then $Q$ has the
reduced generalized c.p.s. structure $\Phi_Q$ of $\Phi$ via
$(N,\nu_\pi N)$ and $\Phi_Q$ is the non degenerate, generalized
c.p.s. structure associated to the Hitchin pair $(\varpi_Q,A_Q)$,
where $\varpi_Q$ is the reduction of $\varpi$ and $A_Q$ is the
projection of $A_N=A|_{TN}$.
\end{corol}
\begin{proof} Under the hypotheses, $\nu_\pi
N=(TN)^{\perp_\varpi}$ and the existence of the reduced symplectic
form $\varpi_Q$ is well known (e.g., \cite{V-carte}, p. 103). Then,
the assertion of the present corollary clearly follows from Theorem
\ref{a3-areducere}. We may also notice that it is easy to justify
the assertion of the corollary straightforwardly. Indeed all we
still need is the fact that $A_N$ sends a $\mathcal{C}(N)$-foliated
vector field $X\in\chi^1(N)$ to a $\mathcal{C}(N)$-foliated vector
field $AX$. In view of the definition of $\mathcal{C}(N)$, this is
equivalent with
$$\varpi([Y,AX],X')=0, \hspace{5mm} \forall X,X'\in\chi^1(N),
\;\forall Y\in\Gamma(TN\cap T^{\perp_\varpi})$$ where $X,X'$ are
$\mathcal{C}(N)$-foliated vector fields, which follows from
$$d\varpi_A(Y,X,X')=0,\;d\varpi(Y,AX,X')=0.$$
\end{proof}

Finally, the following result is a straightforward consequence of
Theorem \ref{a3-areducere}, and may be seen as a Marsden-Weinstein
reduction theorem for generalized c.p.s. manifolds.
\begin{theorem}\label{thMW} Let $(M,\Phi)$ be a generalized c.p.s.
manifold with the Poisson structure $\pi$ and the tensor fields
$A,\sigma$. Assume that one has a $\pi$-Hamiltonian action of the
Lie group $G$ on $M$ with an equivariant momentum map
$J:M\rightarrow\mathcal{G}^*$ ($ \mathcal{G}$ is the Lie algebra of
$G$) such that $J_*\circ A=J_*$. Let $\gamma\in\mathcal{G}^*$ be a
common regular value of all the restrictions of $J$ to the
symplectic leaves of $\pi$ with isotropy group $G_\gamma$. Assume
that the level set $M_\gamma=J^{-1}(\gamma)$ is non void and the
foliation of $M_\gamma$ by the orbits of $G_\gamma$ is by the fibers
of a submersion $s:M\rightarrow Q$. Then the structure $\Phi$
reduces to a generalized c.p.s. structure $\Phi_Q$ of
$Q$.\end{theorem}
\begin{proof} For all the notions involved in Theorem \ref{thMW}
and for the existence of the reduced Poisson structure $\pi_Q$ we
refer the reader to \cite{V-carte}, pp. 110-113. The hypothesis
$J_*\circ A=J_*$ ensures that $M_\gamma$ is $A$-invariant, and the
conclusion follows from Theorem \ref{a3-areducere}.\end{proof}
\begin{rem}\label{obslaMW} {\rm Notice that in Theorem \ref{thMW}
we didn't have to ask the action to be by generalized c.p.s.
mappings, which would have been more restrictive. On the other hand,
if $\pi=0$ the action of $G$ must be trivial, $Q=M$, and we do not
get a true reduction.}
\end{rem}
\section{Generalized c.p.s. submanifolds}
In this section we discuss our second subject, submanifolds. The
naive definition of a generalized c.p.s. submanifold $N$ of a
generalized c.p.s. manifold $(M,\Phi)$ would be by asking the
immersion $\iota:N\hookrightarrow M$ to be a generalized c.p.s.
morphism. Like in Poisson geometry, this condition is very
restrictive because it asks $N$ to be a Poisson submanifold of $M$,
hence, a union of symplectic leaves of the Poisson structure $\pi$
of $\Phi$. The same situation appears if we try to get the
submanifold structure by reducing $\Phi$ via $N$ with control
subbundle $E=0$.

The good notion of a submanifold of a Poisson manifold, which gets
an induced Poisson structure, is that of a Poisson-Dirac submanifold
\cite{CF}. The submanifolds of a generalized, complex manifold with
an induced generalized, complex structure were defined by
Ben-Bassat-Boyarchenko \cite{BB} and, in this section, we discuss
the meaning of the Ben-Bassat-Boyarchenko definition in classical
terms. (A different notion of submanifold, which does not require an
induced structure was studied in \cite{Galt}.)

We begin by recalling the notion of a Poisson-Dirac submanifold. If
$f:N\rightarrow M$ is a differentiable mapping and $L$ is a Dirac
structure on $M$, we obtain a field $f^*(L)$ of maximal isotropic
subspaces of the fibers of $T^{big}N$ by putting
\begin{equation}\label{strDinapoi} f^*(L)_x = \{(X,f^*\alpha)\,/\,
X\in T_xN,\alpha\in T^*_xM,(f_*X,\alpha)\in L_{f(x)}\}\;\;(x\in
N)\end{equation} (e.g., \cite{BR}). The field (\ref{strDinapoi}) may
not be differentiable; if it is, $f$ is called a {\it backward Dirac
map}. If $f$ is the embedding $\iota:N\hookrightarrow (M,L)$ of a
submanifold, and if $L_N=\iota^*(L)$ is differentiable, $L_N$ must
be integrable \cite{C}, $N$ is called a {\it proper submanifold},
and $L_N$ is the {\it induced Dirac structure}.

Particularly, since a Poisson structure $\pi$ may be seen as the
Dirac structure $\{(\sharp_\pi\alpha,\alpha)\,/\,\alpha\in T^*M\}$,
one defines \cite{CF}
\begin{defin}\label{defsbvarPD} {\rm A proper submanifold
$\iota:N\hookrightarrow M$ of a Poisson manifold $(M,\pi)$ such that
the induced Dirac structure is Poisson is called a {\it
Poisson-Dirac submanifold}.}\end{defin}

It was shown in \cite{CF} that the proper submanifold $\iota:N
\hookrightarrow (M,\pi)$ is Poisson-Dirac iff
\begin{equation}\label{cond1dePD} TN\cap\sharp_\pi(ann\,TN)=
TN\cap\nu_\pi N=0.
\end{equation} An equivalent characterization is
obtained by taking the annihilator of (\ref{cond1dePD}), which
yields
\begin{equation}\label{cond2dePD}
ann(\nu_\pi N) + ann\,TN = T^*_NM.
\end{equation}

In \cite{BB} one uses a similar procedure for a definition of a
notion of {\it generalized, complex submanifold} and, in the mean
time, we refer to generalized, complex structures only. Let $\iota:N
\hookrightarrow (M,\Phi)$ be a submanifold of a generalized, almost
complex manifold and let $L\subseteq T^{big}_cM$ be the
$i$-eigenbundle of $\Phi$. Then, $\iota^*(L)$ may be constructed
like in the real case and, if $\iota^*(L)$ is differentiable, we
will say that the submanifold is {\it proper}. If $\Phi$ is
integrable and $N$ is proper, $\iota^*(L)$ is closed by Courant
brackets (like in the real case \cite{C}). However, we may have
$\iota^*L\cap\overline{\iota^*L}\neq0$, and $\iota^*L$ may not be a
generalized, complex structure on $N$. The definition of \cite{BB}
is
\begin{defin}\label{sbvarcomplexegen} {\rm A submanifold $\iota:N
\hookrightarrow M$ of a generalized, almost complex manifold
$(M,\Phi)$ is a {\it generalized, almost complex submanifold} if $N$
is proper and $\iota^*L$ is a generalized, almost complex structure
on $N$, called the {\it induced structure}.}
\end{defin}

The following theorem expresses the conditions of Definition
\ref{sbvarcomplexegen} in classical terms.
\begin{theorem}\label{condptBBcupi} Let $\Phi$ be a
generalized, almost complex structure of matrix form
{\rm(\ref{matriceaPhi})} on $M$ and let $N$ be a submanifold of
$(M,\Phi)$. Then $N$ is a generalized, almost complex submanifold
iff it satisfies the following three conditions: i) $N$ is a
Poisson-Dirac submanifold of $(M,\pi)$, ii) $A(TN)\subseteq
TN+im\,\sharp_\pi=TN\oplus\sharp_\pi(ann\,TN)$, iii) $pr_{TN}\circ
A$, where $pr_{TN}$ is the natural projection of the direct sum of
ii) onto its first term, is differentiable.
\end{theorem}
\begin{proof} The equality included in condition ii) of the theorem
is an immediate consequence of (\ref{cond1dePD}), (\ref{cond2dePD}),
which hold if condition i) holds. Now, let us prove the necessity of
i). The $i$-eigenbundle of $\Phi$, which is the image of
$(1/2)(Id-i\Phi)$, is given by
\begin{equation}\label{eqluiL} L=\{(X-i(AX+\sharp_\pi\xi),
\xi-i(\flat_\sigma X-\xi\circ A))\,/\,X\in TM,\xi\in T^*M\}.
\end{equation} Denote by $\iota$ the immersion of $N$ in $M$.
Using (\ref{eqluiL}) and the natural identification between $T^*N$
and $T_N^*M/ann\,TN$, which represents the covectors of $N$ as
equivalence classes $[\xi]_{ann\,TN}$, we see that the pullback of
$L$ to $N$ is
\begin{equation}\label{eqluiiL} \iota^*L=\{(X-i(AX+\sharp_\pi\xi),
[\xi-i(\flat_\sigma X-\xi\circ A)]_{ann\,TN}\end{equation}$$/\,X\in
TN,AX+\sharp_\pi \xi\in TN\}.$$

On the other hand, if $N$ is a generalized, almost complex
submanifold of $(M,\Phi)$, $\iota^*L$ defines a generalized, almost
complex structure $\iota^*\Phi$ on $N$ and must be of the form
\begin{equation}\label{eqluiiL2} \iota^*L=\{(Y-i(A'Y+\sharp_{\pi'}
[\eta]_{ann\,TN}), [\eta]_{ann\,TN}-i(\flat_{\sigma'}
Y-[\eta]_{ann\,TN}\circ A'))\},
\end{equation} where $A',\pi',\sigma'$ are the elements of the
matrix representation of $\iota^*\Phi$, and $Y\in TN,\eta\in
T^*_NM$. Thus, every pair of the form (\ref{eqluiiL}) is
identifiable with a pair of the form (\ref{eqluiiL2}), and, since
the real part of the equal, vector and covector, components of the
two pairs must be the same, we must have $X=Y$ and $\xi\sim\eta$
modulo $ann\,TN$. The case $X=Y=0$ shows that any $\eta\in T^*_NM$
is equivalent modulo $ann\,TN$ with some $\xi\in ann(\nu_\pi N)$,
i.e., condition (\ref{cond2dePD}) must hold and $N$ is a
Poisson-Dirac submanifold of $(M,\pi)$ with the induced Poisson
structure $\pi'$.

Now, for $Y=X\in TN,\eta\in ann\,TN$,  (\ref{eqluiiL2}) is a pair of
the form
\begin{equation}\label{pereche1} (X-iA'X,-i\flat_{\sigma'}X)
\end{equation} and the
corresponding pair (\ref{eqluiiL}) must be of the form
\begin{equation}\label{pereche2}
(X-i(AX+\sharp_\pi\xi),-i[\flat_\sigma X-\xi\circ
A]_{ann\,TN}),\end{equation} where $\xi\in ann\,TN$ and
$AX+\sharp_\pi\xi\in TN$. The equality of the pairs
(\ref{pereche1}), (\ref{pereche2}) yields
\begin{equation}\label{eqluiA'} A'X=AX+\sharp_\pi\xi,
\end{equation} whence
\begin{equation}\label{condquasitare} A(TN)\subseteq
TN\oplus\sharp_\pi(ann\,TN), \end{equation} which is condition ii)
of the theorem. Furthermore, (\ref{eqluiA'}) implies
\begin{equation}\label{exprluiA'} A'=pr_{TN}\circ A,\end{equation} where
the projection is defined by the decomposition
(\ref{condquasitare}), hence, condition iii) also holds.

For another expression of $A'$ and in order to compute the $2$-form
$\sigma'$ we denote by $\alpha_X\in ann\,TN$ a $1$-form such that
$\sharp_\pi\alpha_X=pr_{\sharp_\pi ann\,TN}AX$. The form $\alpha_X$
is defined up to the addition of a term $\gamma \in
ker\,\sharp_\pi$, i.e., the equivalence class
$[\alpha_X]_{(ann\,TN)\cap(ker\,\sharp_\pi)}$ is well defined. Thus,
for a differentiable vector field $X\in\chi^1(N)$, the
differentiability of $\alpha_X$ is not ensured and may be assumed if
we assume $dim((ann\,TN)\cap(ker\,\sharp_\pi))=const.$ or,
equivalently, $dim(TN + im\,\sharp_\pi)=const.$ Then, $A'$ is given
by \begin{equation}\label{A'cualpha} A'X=AX-\sharp_\pi(\alpha_X).
\end{equation}

Before going on, we notice the following simple result
\begin{lemma}\label{auxiliardem} If conditions i), ii) of
Theorem {\rm\ref{condptBBcupi}} hold, and if
$\gamma\in(ann\,TN)\cap(ker\,\sharp_\pi)$ then $\gamma\circ A\in
ann\,TN$. \end{lemma} \begin{proof} $\forall Z\in TN$ we have
$$\gamma(AZ)=\gamma(\sharp_\pi\alpha_Z)=-\alpha_Z(\sharp_\pi\gamma)=0.$$
\end{proof}

Back to the proof of Theorem \ref{condptBBcupi},  the definition of
$\alpha_X$ shows that $\xi$ of (\ref{eqluiA'}) is of the form
$\xi=-\alpha_X+\gamma$ with
$\gamma\in(ann\,TN)\cap(ker\,\sharp_\pi)$. Accordingly, Lemma
\ref{auxiliardem} implies
$$[\flat_\sigma X-\xi\circ A]_{ann\,TN} = [\flat_\sigma X-\alpha_X\circ
A]_{ann\,TN},$$ and the equality of the pairs (\ref{pereche1}),
(\ref{pereche2}) yields
\begin{equation}\label{bemol'}\flat_{\sigma'}X= [\flat_\sigma X+\alpha_X\circ
A]_{ann\,TN},\end{equation} equivalently,
\begin{equation}\label{exprluisigma'} \begin{array}{l}\sigma'(X,Y) =
(\iota^*\sigma)(X,Y)+\alpha_{X}(AY)\vspace{2mm}\\
=(\iota^*\sigma)(X,Y)-\pi(\alpha_X,\alpha_Y)\hspace{5mm}(X,Y\in
TN).\end{array}\end{equation} Notice that, although $\alpha_X$ is
not uniquely defined, the result of (\ref{exprluisigma'}) is well
defined in view of Lemma \ref{auxiliardem}.

Thus, we proved that a generalized, complex submanifold satisfies
conditions i), ii), iii) and we computed the classical tensor fields
of the induced structure.

For the converse result we first check that i), ii), iii) imply
$\iota^*L\cap\overline{\iota^*L}=0$.

Notice that, if condition i) holds, we have (\ref{cond2dePD}) and
condition ii) is equivalent with (\ref{condquasitare}).

Then, a pair of the form (\ref{eqluiiL}) belongs to
$\iota^*L\cap\overline{\iota^*L}$ iff it also has the form
$$(X+i(AX+\sharp_\pi\zeta),[\zeta+i(\flat_\sigma X-\zeta\circ
A)]_{ann\,TN})$$ with the same vector $X$, i.e., $\exists\zeta\in
T^*_NM$ such that
\begin{equation}\label{eqptAsigma}
\begin{array}{c}
\zeta=\xi+\alpha\hspace{2mm}(\alpha\in ann\,TN),\hspace{2mm}
2AX=-\sharp_\pi(\xi+\zeta),\vspace{2mm}\\ (\xi+\zeta)\circ A -
2\flat_\sigma X\in ann\,TN.\end{array}\end{equation}

Since these conditions imply
$$2(AX+\sharp_\pi\xi)=-\sharp_\pi\alpha,$$ in view of
(\ref{cond1dePD}), we have
\begin{equation}\label{condptAXxi}
\sharp_\pi\alpha=0, \;AX+\sharp_\pi\xi=0,\end{equation} and, if we
apply $A$ to the second condition (\ref{condptAXxi}) and use
(\ref{2-reddezv}), we get
\begin{equation}\label{formaluiX}
X=-\sharp_\pi\flat_\sigma X+\sharp_\pi(\xi\circ A).\end{equation}
Furthermore, using (\ref{eqptAsigma}), (\ref{condptAXxi}) and Lemma
\ref{auxiliardem}, we get
\begin{equation}\label{nrnou}\flat_\sigma X-\xi\circ A\in
ann\,TN.\end{equation} Accordingly, (\ref{cond1dePD}) and
(\ref{formaluiX}) show that $X=0$ and, then, (\ref{condptAXxi}) and
(\ref{nrnou}) show that the pair we study must be of the form
$(0,[\xi]_{ann\,TN})$.

But, $\xi$ is not arbitrary either. Modulo $X=0$,
(\ref{condptAXxi}), (\ref{formaluiX}) and (\ref{nrnou}) imply
$$\xi\in ker\,\sharp_\pi,\;\xi\circ
A\in(ann\,TN)\cap(ker\,\sharp_\pi),$$ and by Lemma \ref{auxiliardem}
$\xi\circ A\in ann\,TN$. Thus, composing again by $A$, we get
$$\xi\circ A^2=-\xi-\flat_\sigma\circ\sharp_\pi\xi=-\xi\in
ann\,TN,$$ and the considered pair is just $(0,0)$. In other words,
we showed that $\iota^*L\cap\overline{\iota^*L}=0$.

The previous conclusion means that, $\forall x\in N$, $\iota^*L_x$
defines a generalized complex structure of $T_x^{big}N$, therefore,
$\iota^*L_x$ must be of the form (\ref{eqluiiL2}). Then, $\pi'$ of
(\ref{eqluiiL2}) is induced by $\pi$ and it is differentiable
because $N$ is a Poisson-Dirac submanifold. Furthermore, $A'$ is
differentiable because of condition iii) and $\sigma'$ is
differentiable because its values depend only on
$\sharp_\pi\alpha_X,\sharp_\pi\alpha_Y$ and the definition of
$\alpha_X,\alpha_Y$ shows that the previous vector field are
differentiable if $A'$ is differentiable. This justifies the fact
that $N$ is proper in $(M,\Phi)$, which finishes the proof of the
theorem. \end{proof}

In view of Theorem \ref{condptBBcupi}, we propose the following
general terminology.
\begin{defin}\label{sbvarcpt} {\rm A submanifold $\iota:N
\hookrightarrow M$ of a generalized, almost c.p.s. manifold
$(M,\Phi)$ is a {\it quasi-invariant submanifold} if i) it is a
$\pi$-Poisson-Dirac submanifold, ii) $A(TN)\subseteq
TN+im\,\sharp_\pi = TN\oplus\sharp_\pi(ann\,TN)$, iii) $A'=
pr_{TN}\circ A$ is differentiable.}\end{defin}

We notice that quasi-invariance is preserved by a gauge equivalence.
Theorem \ref{condptBBcupi} tells that ``quasi-invariant submanifold"
and ``generalized, almost complex submanifold in the sense of
\cite{BB}" are synonymous terms. In the general c.p.s. case the
structure $\Phi$ induces a differentiable, $g$-skew-symmetric
endomorphism $\Phi_N$ of the bundle $T^{big}N$ defined by the matrix
\begin{equation}\label{PhiN}
\Phi_N = \left( \begin{array}{cc}A'&\sharp_{\pi'}\vspace{2mm}\\
\flat_{\sigma'}&-\,^t\hspace{-1pt}A'\end{array}\right),\end{equation}
where $\pi'$ is induced by $\pi$ and $A',\sigma'$ are given by the
formulas (\ref{A'cualpha}), (\ref{exprluisigma'}).
\begin{theorem}\label{strPhiN} If $\iota:N \hookrightarrow
M$ is a quasi-invariant submanifold of a generalized, almost c.p.s.
manifold $(M,\Phi)$, the induced structure $\Phi_N$ is also almost
c.p.s. In the complex and paracomplex cases, if $\Phi$ is
integrable, $\Phi_N$ is integrable too.
\end{theorem}
\begin{proof} For the algebraic part of the proposition it
suffices to work at a fixed point $x\in N$. It is known from
\cite{CF} that, since $N$ is a Poisson-Dirac submanifold of
$(M,\pi)$, there exists a normal space $\nu_xN$
($T_xM=T_xN\oplus\nu_xN$) such that
\begin{equation}\label{Pdescompus} \pi=\pi_{\nu_xN}+\pi_{T_xN},\hspace{1cm}
\pi_{\nu_xN}\in\wedge^2\nu_xN,\,\pi_{T_xN}\in\wedge^2T_xN,\end{equation}
and the induced Poisson structure is $\pi'_x=\pi_{T_xN}$. The
compatibility of $(\pi',A')$ is a straightforward consequence of the
previous remark on $\pi'$, of the compatibility of $(\pi,A)$ and of
formula (\ref{exprluiA'}).

In order to check the compatibility of $(\sigma',A')$ we use formula
(\ref{exprluisigma'}) and, for $X,Y\in T_xN$, we get
$$\sigma'(X,A'Y)= \sigma(X,A'Y)+\alpha_X(AA'Y) = \sigma(X,AY) -
\sigma(X,\sharp_\pi\alpha_Y)$$ $$+\alpha_X(A^2Y) -
\alpha_X(A\sharp_\pi\alpha_Y) \stackrel{(\ref{2-reddezv})}{=}
\sigma(X,AY) - \sigma(X,\sharp_\pi\alpha_Y) +
\sigma(Y,\sharp_\pi\alpha_X)$$ $$ + \pi(\alpha_X\circ A,\alpha_Y).$$
If we change the role of $X,Y$ we see that
$\sigma'(X,A'Y)=\sigma'(A'X,Y)$ as required.

The first condition (\ref{2-reddezv}) for the induced structure is
checked as follows. We have
$$\sharp_{\pi'}\flat_{\sigma'}X= \sharp_{\pi'}[\flat_\sigma
X+\alpha_X\circ A]_{ann\,TN}\hspace{5mm}(X\in TN),$$ and the right
hand side of this equality is computable by a representative form of
the equivalence class sent by $\sharp_\pi$ in $TN$. By
(\ref{2-reddezv}) for the original structure we have
$$\sharp_{\pi}(\flat_\sigma X+\alpha_X\circ A) = \epsilon X-A^2X +
\sharp_\pi(\alpha_X\circ A) $$ $$=\epsilon X-AA'X = \epsilon X -
A'^2X - \sharp_\pi\alpha_{A'X}.$$ Hence,
$$\sharp_{\pi}(\flat_\sigma X+\alpha_X\circ A
+ \alpha_{A'X})\in TN$$ and, since $\alpha_{A'X}\in ann\,TN$, we
deduce that
$$\sharp_{\pi'}\flat_{\sigma'}X =
\epsilon X - A'^2X,$$ which is the required property.

Now, using (\ref{exprluiA'}), (\ref{exprluisigma'}), we get
\begin{equation}\label{nouPhi'} \Phi_N(X,[\xi]_{ann\,TN})=
(AX-\sharp_\pi\alpha_X+\sharp_\pi\tilde\xi, [\flat_\sigma
X+\alpha_X\circ A-\tilde\xi\circ A]_{ann\,TN}),\end{equation} where
$\tilde\xi$ is determined by the decomposition
$\xi=\tilde\xi+\xi_0$, $\tilde\xi\in ann(\nu_\pi N),\xi_0\in
ann\,TN$. Formula (\ref{nouPhi'}) allows us to write down the
general expression of an element of the $\pm i$ or
$\pm1$-eigenbundles, of $\Phi_N$, respectively, similar to the pairs
(\ref{eqluiiL2}) (which was the case of the eigenvalue $i$ in the
complex situation). The results show that these elements also have a
corresponding expression of the form (\ref{eqluiiL}), where, if the
(\ref{eqluiiL2})-like formula is defined by the pair
$(Y=X,[\eta]_{ann\,TN}=[\xi]_{ann\,TN})$, the corresponding
(\ref{eqluiiL})-like formula is defined by the pair
$(X,\tilde\xi-\alpha_X)$.

The conclusion is that the Dirac eigenbundles of the structure
$\Phi_N$ are the $\iota^*$-pullbacks of the Dirac bundles of the
structure $\Phi$. Obviously, under the hypotheses of the theorem,
these pullbacks are differentiable, therefore, if the Dirac
eigenbundles of $\Phi$ are closed by Courant brackets the same holds
for the Dirac eigenbundles of $\Phi_N$. Thus, the assertion about
the integrability of the induced structure follows from Proposition
\ref{integrprinDirac}.
\end{proof}

The proof of the integrability part of Theorem \ref{strPhiN} does
not hold in the generalized, almost subtangent case; the
$0$-eigenbundle of the induced structure $\Phi_N$ is again closed by
brackets, if that of $\Phi$ is (same argument as in the proof
above), but this is not enough for the integrability of $\Phi_N$. If
$\Phi$ is a non degenerate, generalized c.p.s. structure (the
subtangent case included), we can justify the integrability of the
induced structure as follows. The structure $\Phi$ corresponds by
(\ref{corespondenta}) to a Hitchin pair $(\varpi,A)$, and the
submanifold $N$ must be a symplectic submanifold of $(M,\varpi)$ (a
Poisson-Dirac submanifold of a symplectic manifold is a symplectic
submanifold \cite{CF}). Then, (\ref{exprluiA'}) yields
$\varpi'_{A'}=\iota^*\varpi_A$, which is closed. Therefore, the
induced structure $\Phi_N$ is non degenerate, it corresponds to a
Hitchin pair, and it is integrable.

In Section 2, we have defined the notion of an invariant
submanifold, which was a submanifold that is invariant by $A$. Of
course, a Poisson-Dirac, invariant submanifold is quasi-invariant
and, by (\ref{exprluiA'}), (\ref{exprluisigma'}), it has the induced
generalized, almost c.p.s. structure defined by the induced Poisson
structure $\pi'$ and by
\begin{equation}\label{strindinvar}
A'=A|_{TN},\;\sigma'=\iota^*\sigma.\end{equation} In all three
c.p.s. cases, if the structure $\Phi$ of $M$ is integrable, and if
$N$ is a Poisson-Dirac, invariant submanifold, the induced structure
(\ref{strindinvar}) is integrable too. Indeed, the only
integrability condition of Theorem \ref{thCrainic} which is a bit
less obvious is the annulation of the Schouten concomitant
$C_{(\pi',A')}$. But, it is easy to see that $\forall X\in TN$ one
has
\begin{equation}\label{eqfinal2}
<C_{(\pi',A')}([\alpha]_{ann\,TN},[\beta]_{ann\,TN}),X> =
<C_{(\pi,A)}(\tilde\alpha,\tilde\beta),X>,\end{equation} where
$\tilde\alpha,\tilde\beta$ are equivalent modulo $ann\,TN$ with
$\alpha,\beta$ and $\sharp_\pi\tilde\alpha,\sharp_\pi\tilde\beta\in
TN$. Hence, $C_{(\pi,A)}$ implies $C_{(\pi',A')}=0$.

In the terminology of \cite{BB} the invariance of a submanifold $N$
means that the submanifold satisfies the {\it graph condition}. In
\cite{BB} the authors also define a much stronger property called
the {\it split property}. In ``classical terms" a submanifold
$\iota:N\hookrightarrow(M,\Phi)$ is a {\it split submanifold} if it
is Poisson-Dirac, invariant and has an $A$-invariant normal bundle
$\nu N$ $(T_NM=TN\oplus\nu N)$ which is $\sigma$-orthogonal to $TN$.
Like any invariant submanifold, a split submanifold has the induced
generalized structure $\Phi_N$ defined by (\ref{strindinvar}). But,
the normal bundle $\nu N$ also has an induced generalized structure
$\Phi_\nu$, and $\Phi=\Phi_N\oplus\Phi_\nu$.

Finally, let us also make the following observation. The definitions
of quasi-invariance and invariance may also be used for a
submanifold $N$ of a Poisson-Nijenhuis manifold $(M,\pi,A)$
($(\pi,A)$ is a Poisson-Nijenhuis structure). A Poisson-Dirac,
invariant submanifold inherits an induced structure
$(\pi',A'=A|_{TN})$ and formula (\ref{eqfinal2}) shows that the
induced structure is a Poisson-Nijenhuis structure too. Moreover, it
is easy to see that the Poisson-Nijenhuis hierarchy
$(\pi'_k,A^{'p})$ ($k,p=1,2,...$,
$\sharp_{\pi'_k}=A^{'k}\circ\sharp_{\pi'}$) of the induced structure
$(\pi',A')$ is induced by the corresponding Poisson-Nijenhuis
hierarchy of $(\pi,A)$. For these reasons it is natural to attribute
the name of {\it Poisson-Nijenhuis submanifold} to an invariant,
Poisson-Dirac submanifold $N$ of a Poisson-Nijenhuis manifold
$(M,\pi,A)$. If $N$ is a quasi-invariant submanifold of $(M,\pi,A)$,
$N$ has the induced Poisson structure $\pi'$ and a compatible,
tensor field $A'$ defined by (\ref{exprluiA'}) but $A'$ may not have
a vanishing Nijenhuis tensor.
\end{document}